# Resolution-independent generative models based on operator learning for physics-constrained Bayesian inverse problems

Xinchao Jiang[1, a], Xin Wang[a], Ziming Wen[a], Hu Wang[*, a, b]

*a. State Key Laboratory of Advanced Design and Manufacturing Technology for Vehicle, Hunan University, Changsha, 410082, P.R. China*

*b. Beijing Institute of Technology Shenzhen Automotive Research Institute, Shenzhen, 518000, P.R. China*

## Abstract

The Bayesian inference approach is widely used to tackle inverse problems due to its versatile and natural ability to handle ill-posedness. However, it often faces challenges when dealing with situations involving continuous fields or large-resolution discrete representations (high-dimensional). Moreover, the prior distribution of unknown parameters is commonly difficult to be determined. In this study, an Operator Learning-based Generative Adversarial Network (OL-GAN) is proposed and integrated into the Bayesian inference framework to handle these issues. Unlike most Bayesian approaches, the distinctive characteristic of the proposed method is to learn the joint distribution of parameters and responses. By leveraging the trained generative model, the posteriors of the unknown parameters can theoretically be approximated by any sampling algorithm (e.g., Markov Chain Monte Carlo, MCMC) in a low-dimensional latent space shared by the components of the joint distribution. The latent space is typically a simple and easy-to-sample distribution (e.g., Gaussian, uniform), which significantly reduces the computational cost associated with the Bayesian inference while avoiding prior selection concerns. Furthermore, incorporating operator learning enables resolution-independent in the generator. Predictions can be obtained at desired coordinates, and inversions can be performed even if the observation data are misaligned with the training data. Finally, the effectiveness of the proposed method is validated through several numerical experiments.

---

[1] First author. *E-mail address*: jiangxinchao@hnu.edu.cn (X.C. Jiang)

[*] Corresponding author. *E-mail address*: wanghu@hnu.edu.cn (H. Wang)

***Keywords*:** Operator learning; generative adversarial networks; resolution-independent; Bayesian; inverse problems.

## 1. Introduction

Many research scenarios revolve around a similar problem of identifying the input parameters (causes) of physical systems, such as sources, boundary conditions, material parameters, etc. These problems are commonly referred to as parameter identification or inverse problems. Most of these so-called inverse problems are typically ill-posed according to Hadamard's definition [1]. This includes solution existence, solution uniqueness, and instability of the solution procedure. In practice, most difficulties arise from the fact that, when going against causality, forward stability turns into inverse instability since observations (consequences) are inevitably corrupted by noise [2].

The inverse problem is often well-posed by regularization and reformulated deterministically before it can be solved using optimization or other means [3]. Many approaches have also been proposed within the framework of statistical inference to address this problem [4]. However, this study focuses specially on Bayesian methods [5], which provide a complete probabilistic description of unknown parameters given observed data and obtain the posterior distribution of quantities of interest (QoIs) by merging likelihood and prior information. It provides a robust manner to handle system variabilities and parameter fluctuations. And it is one of the most popular approaches to infer probabilistic descriptions of model parameters. Certainly, classical Bayesian methods represent a major advance in the field of inverse problems. The output from Bayesian inference provides quantitative estimates for the distribution of the identification results, which can be critical in applications involving high-stake decisions. However, certain important issues related to prior distributions, computational cost, likelihood functions, and others have to be solved under the widely used framework first. Specifically, two particular issues are covered in this study.

1) Choosing an appropriate informative prior distribution for the unknown parameters can be challenging, as inappropriate priors may lead to biased posteriors. Typically, the choice of the explicit prior distribution is limited to several typical distributions, which ensures that the posterior can be obtained analytically or numerically. However, it also limits the expressiveness of the posterior.

2) Exploring high-dimensional parameter space using Markov Chain Monte Carlo (MCMC) methods [6, 7] would be computationally prohibitive and difficult to implement.

To solve above mentioned challenges, various approaches have been explored from different aspects. Such as, some studies use spectral methods such as Karhunen-Loeve Expansion (KLE) [8] to reduce high-dimensional parameters into several coefficients and then explore the posterior based on sampling these coefficients [9, 10]; Some attempts rely on reduced models or surrogate models to speed up the forward simulation in computationally intensive inverse problems [11, 12]. In another technique route called Approximate Bayesian Computation (ABC), researchers consider using some distance metric (e.g., Euclidean distance) based on the summary statistics to replace the intractable likelihood function, and typically use a dimension reduction algorithm to project observations into a low-dimensional manifold to accelerate the computation [13-15]; Recently, a novel method called Physics-informed Neural Network (PINN) [16, 17] is attractive. Compared with other popular methods, the distinctive characteristics are to avoid the problem of repeated forward simulations of the inverse problem and to reduce data requirements by introducing physical laws. However, the aforementioned methods represent non-generative models that typically focus on the conditional distribution of causes for given consequences. In recent years, the deep generative model has received more attention, and it has also rapidly developed in solving inverse problems [18-21]. For instance, Generative Adversarial Networks (GANs) [22] are used to learn prior distributions of complex parameters in Bayesian inference [23]. Some studies are centered around learning the joint distribution of causes and consequences, showcasing the flexibility of joint distributions in problem-solving [24-26]. Among them, the work of Patel et al. [24] has shown how to use the approximate joint distribution learned by GANs. This framework solves the issue of prior selection by using GANs and significantly reduces the inference cost by transferring the original problem to a low-dimensional latent space where each component of the joint distribution shares these latent variables. While impressive, the computational physics inverse problem solved in their study is resolution-dependent, meaning the deep learning model input and output dimensions are fixed. In computational physics, the number of pixels is not always suitable as the dimension of the physical field, which would result in redundant dimensions and an unfavorable selection of latent variable dimensions.

Therefore, to address the aforementioned issues in the GAN-based inverse problem, we have introduced a recent line of deep learning called Deep Operator Network (DeepONet) [27]. In practice, the Wasserstein GAN with Gradient Penalty (WGAN-GP) [28] is used. By incorporating the coordinate space into the generative model, the space of interest (parameters, or responses) becomes resolution-independent. The generator can produce responses for any positions within the domain by combining the latent variables and the coordinates. This feature allows the observed data can be misaligned with the training data, which is critical when dealing with sensor limitations such as number and installation position. Additionally, assuming the noise of the observed data is an additive Gaussian, and its statistical properties are also known, therefore the problem of the likelihood function can be solved naturally. Finally, the effectiveness of the proposed method is demonstrated by several numerical examples, where the inverse problem is addressed by using only a few latent variables.

The remainder of this paper is organized as follows. In Section 2, the formulation of the Bayesian inverse problem based on joint distribution is given and the key issue of this framework is mentioned. Thereafter, the details of the proposed method are described in Section 3. Several numerical examples are presented in Section 4. Finally, the core of this study is summarized in Section 5. Upon formal publication of the paper, all code and data associated with this study will be made available to the public at

## 2. Problem formulation

Bayesian inference is a probabilistic framework that determines causes from desired or observed consequences. This framework starts from the observation and reveals that the observations are subject to noise and that it is common to formalize this as an operator equation

$$\hat{\boldsymbol{d}} = \mathcal{G}(\boldsymbol{m}) + \eta, \tag{1}$$

where $\hat{\boldsymbol{d}}$ is the given noisy observation data. $\mathcal{G}(\cdot)$ is the function or operator. $\boldsymbol{m}$ is the parameters. $\eta$ is a mean-zero additive noise whose statistical properties $p_\eta$ are known. The task of the classical Bayesian inverse problem is to recover the parameters $\boldsymbol{m} \in \mathcal{M} \subset \mathbb{R}^{N_\mathcal{M}}$ from measured data $\hat{\boldsymbol{d}} \in \mathcal{D} \subset \mathbb{R}^{N_\mathcal{D}}$, where the domains are compact, and its task can be considered as learning the posterior $p_\mathcal{M}(\boldsymbol{m}|\hat{\boldsymbol{d}})$. However, the focus of this study is on a more recent formulation [24] that differs significantly from classical Bayesian methods by using the joint distribution of the parameters and responses $\boldsymbol{a} = [\boldsymbol{m}, \boldsymbol{d}] \in \mathcal{A} \equiv \mathcal{M} \times \mathcal{D} \subset \mathbb{R}^{N_\mathcal{M} + N_\mathcal{D}}$. The approximated joint distribution is constructed based on the deep generative models. Once the joint distribution is obtained, the corresponding conditional distribution and marginal distribution can be easily obtained. Therefore, this formulation can flexibly deal with different problems in one model, including forward problems, inverse problems, and mixed problems. While in this study, the focus is placed on its application to the inverse problem, which can be described by Bayes' rule:

$$p_\mathcal{A}^{\text{post}}(\boldsymbol{a}|\hat{\boldsymbol{d}}) = \frac{1}{\mathbb{Z}} p^{\text{like}}(\hat{\boldsymbol{d}}|\boldsymbol{a}) p_\mathcal{A}^{\text{prior}}(\boldsymbol{a}), \tag{2}$$

where $p^{\text{like}}$ is the likelihood of $\hat{\boldsymbol{d}}$ given $\boldsymbol{a}$. $\mathbb{Z}$ is the normalization constant which is also referred to as model evidence. Consider the mathematical-physical problem which is well-posed and given the condition of a definite solution, and assume that the *m* and *d* correspond one to one. Consequently, denoting $\hat{d} = t_d(a) + \eta$, where $t_d(\cdot)$ means the operation of extracting *d* from *a*. The equation (2) can be rewritten as:

$$p_{\mathcal{A}}^{post}(a|\hat{d}) = \frac{1}{\mathbb{Z}} p_\eta(\hat{d} - t_d(a)) p_{\mathcal{A}}^{prior}(a). \tag{3}$$

This method does not directly compute the posterior of the unknown parameters but obtains the posterior $p_{\mathcal{M}}^{post}(m|\hat{d})$ in two steps. First, the posterior samples of the joint distribution *a* are obtained. Then the marginal distribution of the parameters *m* is extracted as

$$m = t_m(a), \tag{4}$$

where $t_m(\cdot)$ is the operation of extracting *m* from *a*. The key issue of computing the posterior density in equation (3) requires the construction of the joint prior distribution $p_{\mathcal{A}}^{prior}(a)$. The approach to constructing such a prior distribution will be covered in detail in Section 3 based on an available data set $\mathcal{S} = \{a_1, a_2 \ldots, a_s\}$. The statistics (e.g., mean, standard deviation) of the posterior distribution of *a* can be approximated by MCMC sampling.

## 3. Methodology

This section presents the details of the proposed method. Figure 1 is the flowchart that summarizes the core of the methodology. The main building blocks involved are described separately.

### 3.1 *Generative adversarial networks*

GANs have been proven to be successful in learning the probability distribution from data. GANs deal with this problem by defining a 2-player zero-sum game between the generator $G(\cdot)$ and the discriminator $D(\cdot)$, which are two Deep Neural Networks (DNNs) parameterized by different parameters. Consider the problem of approximating the distribution $p_{\mathcal{A}}^{prior}$ based on $\mathcal{S}$. The purpose of the generator is to map a latent vector $z \in \mathcal{Z} \subset \mathbb{R}^{N_\mathcal{Z}}$ to $\mathcal{A}$, where $N_\mathcal{Z} \ll N_\mathcal{M} + N_\mathcal{D}$. The input $z$ in the latent space is typically sampled from a simple distribution (e.g., uniform, Gaussian), and the distribution of the push-forward of latent variables through the generator $G(z)$ is denoted as $p_\mathcal{A}^G$. Meanwhile, the discriminator is employed to evaluate the discrepancy between the $p_\mathcal{A}^{\text{prior}}$ and $p_\mathcal{A}^G$ to accomplish such a game, whose input contains samples from both two distributions. The generator and the discriminator are trained in an adversarial manner. Different types of GANs can be trained by appropriate loss functions. In this work, the WGAN-GP is used since it hardly suffers from the problem of mode collapse problem and the training is stable [28]. Specifically, the loss function for the generator and discriminator (critic) of the WGAN-GP is defined as:

$$\mathcal{L}(D, G) = \mathbb{E}_{a \sim p_\mathcal{A}^{\text{prior}}}[D(a)] - \mathbb{E}_{z \sim p_\mathcal{Z}}[D(G(z))] + \lambda \mathbb{E}_{\hat{a} \sim p_{\hat{\mathcal{A}}}}[(\|\nabla_{\hat{a}} D(\hat{a})\|_2 - 1)^2], \qquad (5)$$

where $\lambda$ is the penalty coefficient. $p_{\hat{\mathcal{A}}}$ is the distribution generated by uniform sampling on straight lines between pairs of points from $p_{\mathcal{A}}^{\text{prior}}$ and $p_{\mathcal{A}}^{G}$. The gradient penalty soft constraint is added to the discriminator to ensure the 1-Lipschitz condition. The generator is trained to minimize the equation (5), while the discriminator is trained to maximize it. It can be concluded as a min-max problem to determine the optimal solution:

$$(D^*, G^*) = \inf_G \sup_D \mathcal{L}(D, G). \tag{6}$$

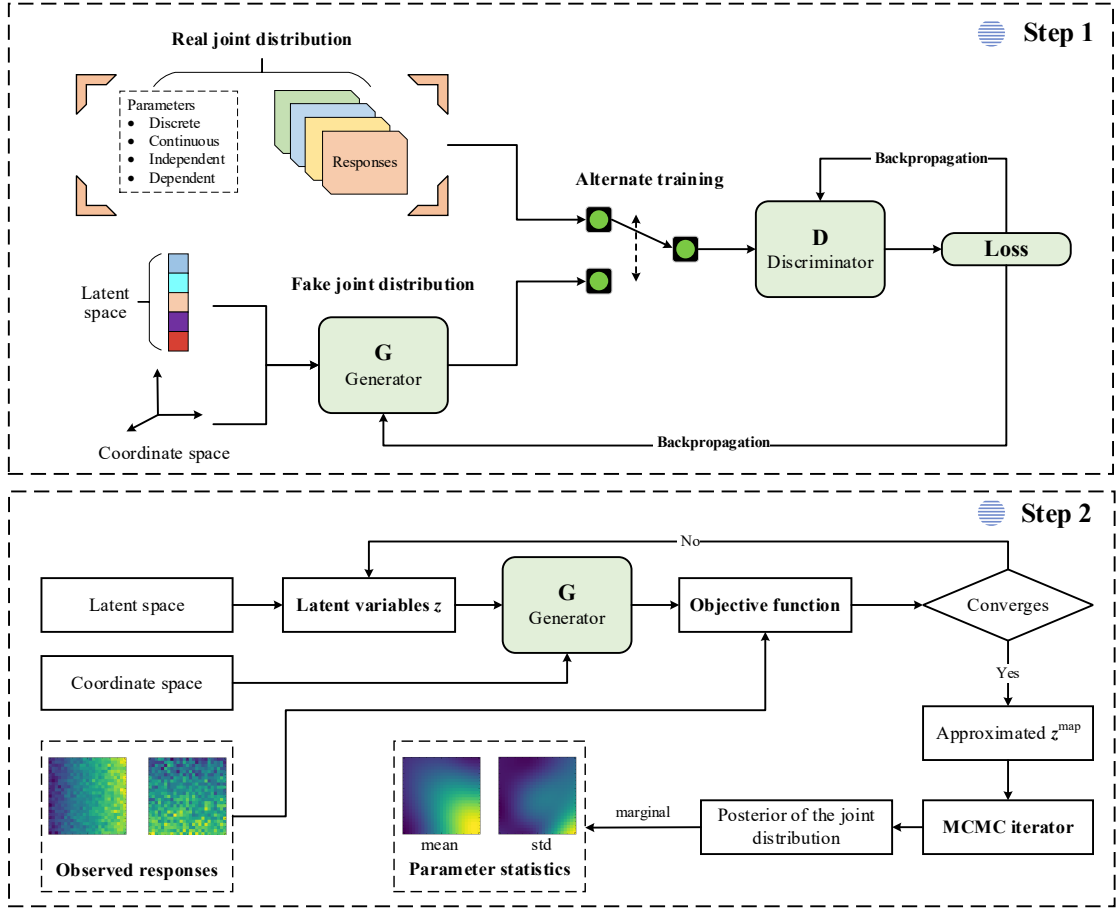

**Figure 1** The framework of the proposed inverse method. In Step 1, a generator is trained based on the GANs framework. In Step 2, with given observations, optimization and MCMC iteration are implemented based on the trained generator in Step 1.

The inner maximization problem corresponds to the discriminator which is trained under the 1-Lipschitz constraint and results in an approximation of the Wasserstein-1 distance [29] between $p_{\mathcal{A}}^{\text{prior}}$ and $p_{\mathcal{A}}^{G}$. The outer minimization problem refers to finding a generator that minimizes the Wasserstein-1 distance. Finally, the push-forward of $p_{\mathcal{Z}}$, $G^*(z)$, can be viewed as the approximation of the prior data distribution.

### 3.2 Deep operator networks

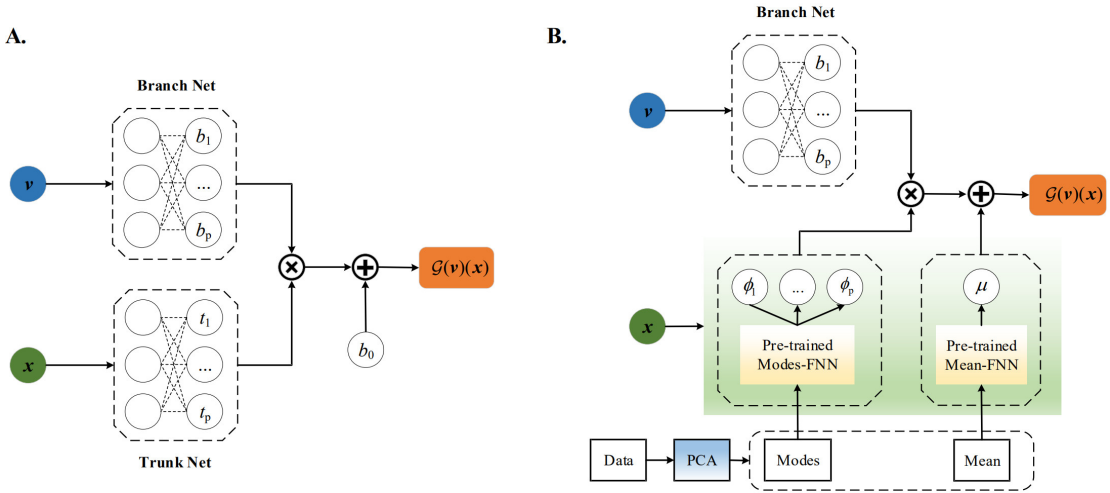

**Figure 2** Illustrations of the architecture of DeepONet (A) and POD-DeepONet (B).

It should be noticed that the operator networks represented by DeepONet have achieved great success in recent years. In this section, a brief overview of DeepONet is presented. Considering a physical system, the aim of DeepONet can be viewed as learning solution operators of parametric PDEs. In this framework, it consists of two subnetworks that have two independent inputs. One is the branch network, which generates the *p*-dimensional vector ***b*** by encoding the discretized input functions or parameters, ***v***, at given sensor positions; The other is the trunk network, which generates the *p*-dimensional vector ***t*** by encoding the coordinates ***x***. Two forms of this framework are mainly introduced in this work. In the vanilla unstacked DeepONet (vanilla refers to the unmodified version), the final output of the operator network is expressed as:

$$\mathcal{G}(\bm{v})(\bm{x}) \approx \sum_{k=1}^{p} b_k(\bm{v}) t_k(\bm{x}) + b_0, \qquad (7)$$

where $b_0 \in \mathbb{R}$ is a bias (trainable). $\{b_1, b_2, \ldots, b_p\}$ are the *p* outputs of the branch network, and $\{t_1, t_2, \ldots, t_p\}$ are the *p* outputs of the trunk network. The trunk network automatically learns the basis of the output function from the data. While in the POD-DeepONet [30], proper orthogonal decomposition (POD) or principal component analysis (PCA) [31] is first imposed on the training data, and the trunk network is replaced by a pre-trained network based on the modes $\phi_k(\bm{x})$ of the data. The bias term is replaced by a pre-trained network of the mean $\mu(\bm{x})$ of the data. Specifically, the output can be described as:

$$\mathcal{G}(\bm{v})(\bm{x}) \approx \sum_{k=1}^{p} b_k(\bm{v}) \phi_k(\bm{x}) + \mu(\bm{x}). \qquad (8)$$

Notably, the parameters of these two pre-trained surrogate models are frozen in the training phase of POD-DeepONet. Moreover, for both types of operator networks, the output $\mathcal{G}(\bm{v})(\bm{x})$ can be viewed as a function of ***x*** conditioning on ***v***.

### 3.3 *Operator Learning-based GANs: OL-GANs*

In general, GANs that are trained based on a discrete representation of data (pixels, sparse sensors) are resolution-dependent. For a fixed dimension generator, the generated data is also based on the pixels or coordinate points of the training set. Therefore, the generator cannot give a direct prediction for a new position, i.e., resolution-dependent. In addition, determining the dimension of latent variables required for high-resolution physical field data is difficult. The redundant dimension will cause the cost of network training and the computational cost of the subsequent Bayesian inference is still expensive. In this study, to address the aforementioned issues, we proposed a novel construction of generators of GANs in some situations motivated by the architecture and theory of DeepONets. This idea has a similar application in [32]. Specifically, a physical system is often described by a Cartesian coordinate system. Thus, generators of the proposed method for continuous physical fields can be described as:

$$G(z)(x) = \sum_{k=1}^{p} b_k(z) t_k(x) + b_0, \tag{9}$$

or POD-based:

$$G(z)(x) = \sum_{k=1}^{p} b_k(z) \phi_k(x) + \mu(x). \tag{10}$$

where the parameters (i.e., discrete input function or features) are replaced by a random vector, $z$, from independent identically Gaussian distributions in the latent space. For clarity, the GANs proposed in this study with such architectures are denoted as OL-GANs which include two parts of inputs, that is the coordinate and the latent vector. The DeepONet typically constructs a discriminative model based on the mean square error (MSE). However, this study constructs a generative model based on Wasserstein-1 distance to approximate $p_{\mathcal{A}}^{\text{prior}}$. For any point $x$ in the domain of $G(z)$, the output $G(z)(x)$ is a real number corresponding to a specific coordinate. Given a set of coordinates **x** concerning **a**, the loss of the OL-GANs is rewritten as:

$$\mathcal{L}(D, G) = \mathbb{E}_{a \sim p_{\mathcal{A}}^{\text{prior}}}[D(a)] - \mathbb{E}_{z \sim p_{\mathcal{Z}}}[D(G(z)(\mathbf{x}))] + \lambda \mathbb{E}_{\hat{a} \sim p_{\hat{\mathcal{A}}}}[(\|\nabla_{\hat{a}} D(\hat{a})\|_2 - 1)^2]. \tag{11}$$

It is worth noting that generator architecture does not have to use this trunk and branch network-based design for components that are not described in coordinate systems. In other words, the input **x** is not necessary in such cases, and the equation (5) still serves as the loss function. The primary advantage of incorporating the coordinate space into the generator lies in rendering it resolution-independent. Generators can learn the solution operators from random variables to continuous physical fields instead of discrete representations of the field. For example, when presented with new sensor observation data, all relevant coordinates would be $\mathbf{x}_{\text{new}}$. The trained OL-GAN can directly generate responses at this new position, $G(z)(\mathbf{x}_{\text{new}})$, and the subsequent inverse procedure can still be performed, without the dilemma that the fake data are misaligned with the observation data. In addition, generative models transfer the original parameter and response space to the latent space, which could greatly reduce the computational cost associated with subsequent inference tasks. According to the reference [33], the surrogate construction approach of DeepONet relies on finding a suitable projection that approximates the operator with good accuracy. We recommend using projection methods such as PCA to estimate the appropriate output dimensions, $p$, for both the branch network and trunk network. The appropriate dimensions of $z$ can be selected in the range of less than $p$. Since the dimensionality of the latent variables and coordinates are low, all the network architectures in this study are based on simple Fully-connected Neural Networks (FNNs). However, other network architectures can also be used depending on the specific problem.

### 3.4 OL-GANs assisted Bayesian inference

Considering an approximated distribution $p_{\mathcal{A}}^{G}$ constructed via the OL-GANs that is weakly equivalent to $p_{\mathcal{A}}^{\text{prior}}$. Assuming the latent variables $z \sim P_{\mathcal{Z}}(z)$ and the push-forward of $z$ at given coordinates $\mathbf{x}$ which is denoted as $G^*(z)(\mathbf{x})$. If $G^*$ is the limit of a class of generators for which satisfy the condition that the Wasserstein-1 distance between training data and generated data is close to zero. Hence, for all continuous bounded functions $f \in \mathcal{C}_b(\mathcal{A})$, there is a relationship between $a$ and $z$ [24]:

$$\mathbb{E}_{a \sim p_{\mathcal{A}}}[f(a)] = \mathbb{E}_{z \sim p_{\mathcal{Z}}}[f(G^*(z)(\mathbf{x}))]. \tag{12}$$

Choosing $f(a) = \frac{1}{\mathbb{Z}} l(a) p_{\eta}(\hat{d} - t_d(a))$, and assuming that $p_{\eta}, l(\cdot) \in \mathcal{C}_b(\mathcal{A})$ too, the equation (12) can be rewritten as:

$$\mathbb{E}_{a \sim p_{\mathcal{A}}}[\frac{1}{\mathbb{Z}} l(a) p_{\eta}(\hat{d} - t(a))] = \mathbb{E}_{z \sim p_{\mathcal{Z}}}[\frac{1}{\mathbb{Z}} l(G^*(z)(\mathbf{x})) p_{\eta}(\hat{d} - t_d(l(G^*(z)(\mathbf{x}))))], \tag{13}$$

$$\mathbb{E}_{a \sim p_{\mathcal{A}}^{\text{post}}}[l(a)] = \mathbb{E}_{z \sim p_{\mathcal{Z}}^{\text{post}}}[l(G^*(z)(\mathbf{x}))], \tag{14}$$

where

$$p_{\mathcal{Z}}^{\text{post}}(z|\hat{d}) \equiv \frac{1}{\mathbb{Z}} p_{\eta}(\hat{d} - t_d(G^*(z)(\mathbf{x}))) p_{\mathcal{Z}}(z). \tag{15}$$

The equation (14) implies that sampling from the posterior distribution of $a$ is equivalent to calculating the push-forward $G^*(z)(\mathbf{x})$ by sampling from the posterior distribution of $z$. That is:

$$a \sim p_{\mathcal{A}}^{\text{post}}(a|\hat{d}) \Leftrightarrow a = G^*(z)(\mathbf{x}), z \sim p_{\mathcal{Z}}^{\text{post}}(z|\hat{d}). \tag{16}$$

Hence, the original Bayesian inference can be implemented in the latent space. This would circumvent the possible "curse of dimensionality" of performing Bayesian inference and the difficulty of a prior selection in the original parameter space. Moreover, it greatly improves the efficiency since the dimension of $z$ is far less than that of $a$. Practically, an efficient approximation for the right-hand side of the equation (14) can be conveniently achieved by Metropolis-Hasting (MH) sampling, which updates based on the logarithmic definition presented in equation (15). We used an adaptive proposal Gaussian distribution to achieve a desired acceptance rate between 0.2 and 0.5 [34]. Consequently, the mean ($\mu_m$) and standard deviation ($\sigma_m$) of the posterior distribution of parameters ($m$) can be described below by choosing $l(a) = t_m(a)$ and $l(a) = t_m(a^2)$:

$$\mu_m = \mathbb{E}_{m \sim p_{\mathcal{M}}^{\text{post}}}[m] = \mathbb{E}_{a \sim p_{\mathcal{A}}^{\text{post}}}[t_m(a)] = \mathbb{E}_{z \sim p_{\mathcal{Z}}^{\text{post}}}[t_m(G^*(z)(\mathbf{x}))], \tag{17}$$

$$\mu_{m^2} = \mathbb{E}_{m \sim p_{\mathcal{M}}^{\text{post}}}[m^2] = \mathbb{E}_{a \sim p_{\mathcal{A}}^{\text{post}}}[t_m(a^2)] = \mathbb{E}_{z \sim p_{\mathcal{Z}}^{\text{post}}}[t_m((G^*(z)(\mathbf{x}))^2)], \tag{18}$$

$$\sigma_m = \sqrt{\mu_{m^2} - \mu_m^2}. \tag{19}$$

Generally, the convergence of the MCMC method is closely tied to the initial point selection, which lacks a definitive approach. According to equation (15), it is possible to determine the maximum a posterior (MAP) estimate first. And then take the $z^{\text{map}}$ as the initial point in the MH sampling procedure. Since the noise has been assumed to be an isotropic Gaussian distribution with zero mean and a variance matrix $\Sigma$. Therefore, the MAP estimate of latent variables, $z^{\text{map}}$, can be obtained by minimizing the objective function through appropriate optimization algorithms. Specifically, the form of the objective function is:

$$\mathcal{L}^{\text{map}} = -\log p_{\mathcal{Z}}^{\text{post}}(z|\hat{d}) \propto \frac{1}{2}\left|\Sigma^{-1/2}(\hat{d} - t_d(G^*(z)(\mathbf{x})))\right|^2 + \frac{1}{2}|z|^2. \tag{20}$$

In practice, $\mathcal{L}^{\text{map}}$ is optimized by Adam [35]. Since the generator is a neural network, it is straightforward to use the mechanism of auto differentiation [36] to obtain the gradient of latent variables $\nabla_z \mathcal{L}^{\text{map}}$. Consequently, the MAP estimation can be obtained rapidly and conveniently.

## 4. Numerical results and discussion

In this section, several numerical examples are presented to demonstrate the effectiveness of the suggested inverse approach. Including both the independent and dependent parameters (e.g., function), these examples show the use of operator learning and coordinate space for different problems. All the code is written in Python, the deep learning part is based on Pytorch [37], and the finite element method used to generate data is based on FEniCS [38]. Some of the hyperparameters do not change in all examples. For instance, the latent variable obeys a standard Gaussian distribution. The gradient penalty coefficient is 10. The number of discriminator iterations per generator iteration is 5. The Kaiming uniform initializer [39] is used to initialize all FNNs used. The Adam is used to optimize the parameters (weights and bias) of generators and discriminators with a fixed learning rate of 1e-4, the training epochs are 20,000, and $\beta_1 = 0, \beta_2 = 0.999$. In OL-GANs, the number of outputs of the trunk net and the branch net, $p$, is 10. The negative slope of the LeakyReLU activation function is 0.01. As for the MAP estimation in the inversion phase, the number of iterations of Adam is 1,000 and the learning rate is 0.05. The number of samples in the MCMC iteration is 10,000, and the burn-in period is 5,000. The settings of all these parameters are empirical parameters obtained based on many trials, and they all work well in the cases involved in this work.

### 4.1 *Case 1: heat source inversion problem*

Starting with a parameter identification problem for a Poisson equation which arises in many fields such as the modeling of subsurface flow and heat conduction. This study focuses on the inverse heat conduction problem where the aim is to identify the unknown parameters that arise in the source function, which is assumed to be a Gaussian function. For a unit square domain $\Omega = [0,1] \times [0,1] \in \mathbb{R}^2$ with a boundary $\partial \Omega = \Gamma_D \cup \Gamma_N$, the governing equation with particular boundary conditions (BCs) can be described as follows:

$$\begin{cases} -\nabla^2 u = f, \text{ in } \Omega, \\ \quad u = 0, \text{ on } \Gamma_D, \\ \nabla u \cdot \mathbf{n} = \sin 5x, \text{ on } \Gamma_N. \end{cases} \quad (21)$$

Where $f = 10\exp(-\frac{1}{2}((x-c_1)^2 + (y-c_2)^2)/c_3^2)$ is the source with variable parameters $c_1, c_2 \in [0.2, 0.8]$ and $c_3 \in [0.05, 0.15]$. $\Gamma_D = \{(0,y) \cup (1,y) \subset \partial\Omega\}$ is the Dirichlet boundary and $\Gamma_N = \{(x,0) \cup (x,1) \subset \partial\Omega\}$ is the Neumann boundary. $\mathbf{n}$ denotes the outward-directed boundary normal. A data set of 1,000 samples of $f$ is generated by Latin Hypercube Sampling (LHS) in the defined uniform distribution domain. The corresponding temperature field is discretized on a $33 \times 33$ grid and solved by FEM. The responses are collected from 81 sparse sensors at the locations by interpolation, which is described in Figure 3. Hence, the data form of the training sample of this case is the joint distribution which is a vector consisting of 3 independent parameters and 81 temperature responses, that is $\boldsymbol{a} = [\boldsymbol{m}, \boldsymbol{d}] = [c_1, c_2, c_3, u_1, \ldots, u_{81}]$. Then, an observation of the temperature field is taken from a test set that is different from any sample in the training set, and it is disturbed by additive Gaussian noise, which is denoted as $\hat{\boldsymbol{d}}$. According to the formulation in Section 2, the task of this case is to infer the 3 parameters concerning $\hat{\boldsymbol{d}}$ by the suggested method, here the true parameters are 0.4489, 0.7340, and 0.1111, respectively. The details of the OL-GANs of this case are shown in Figure 4, the MSE of the pre-trained POD surrogate is less than 1e-8. The generator and the discriminator are optimized by Adam. The batch size is 250. In this case, 3 latent variables are used for the generative models. Two architectures of DeepONets are considered in this study, and the specific computational results for this case are as follows.

The inversion results of 81 sparse uniform sensors with different noise levels are shown in Table 1. The vanilla version means that the trunk net is not POD-based. In this example, the proposed method achieves good results under two different types of trunk nets, the deviation of the posterior mean is small, the maximum relative error is less than 3%, and the truth values are all within the 3-standard deviation (3-std) confidence interval. Figure 5 to Figure 8 are the detail of the posterior distribution of the parameters, where the red dot of the scatter plot represents the position of the posterior mean. The specific observed response is shown in Figure 9. Since we focus on the inversion of the 3 independent parameters, only the relevant information of the posterior response generated by the vanilla OL-GAN is given here for demonstration. When the noise level is 0.01, the accuracy of the two generator architectures is comparable. When the noise is increased to 0.1, the POD-based architecture could achieve relatively better results, but the accuracy improvement is small.

To illustrate that the generator can generate high-quality physical field responses, Figure 10 is drawn. The eigenvalues of the real samples (training set) and fake samples (2000, default) are calculated and normalized by PCA. It can be seen that the explained variance ratio of the fake samples of the two generator types is very similar to that of the real samples, indicating that the generator has learned the underlying distribution of the physical response. It should be noted that the transformation of latent variables into parameters in this example is based on plain FNN, and operator learning can be only applied in physical field responses that can be described by coordinates. Therefore, POD is also only available in networks that correspond to the prediction of the response. In this example, 3 latent variables are used to generate a joint distribution of 84 dependent variables. This is because of the presence of 3 independent parameters and the finding, based on many experiments, that 3 latent variables are sufficient to ensure acceptable accuracy. If more latent variables are used, the process of generating parameters from latent variables is the "parameter dimensionality increase". However, it needs also to be comprehensively considered with the complexity of the responses. Since the dimension of latent variables is much lower than that of the joint distribution, it is worth considering an appropriate increase in the dimension of latent variables to ensure the accuracy of the generated joint distribution of parameters and responses in some situations. In fact, it is possible to obtain feasible results with fewer latent variables, which involves dimensionality reduction of the independent parameters and is discussed in Appendix A.

In addition, the reference [23] uses the ratio of the number of joint distribution variables to the number of latent variables to describe the dimension reduction value. We believe that the number of pixels is not always suitable as the dimension of the physical field in computational physics, and the dimension reduction ratio of this example would be 84/3=28 in this way. On the one hand, in terms of the principal components of the response (linearly), the temperature field in the training data can retain 94.5% of the information with 3 principal components, 99% with 5 components, and 99.9% with 12 components. On the other hand, it is resolution-dependent. High-resolution images/sensors could lead to redundant dimensions. However, the suggested method in this study is resolution-independent, the input of the generator includes the coordinate in addition to the latent variable. A more intuitive explanation is that in the case of fixed latent variables, the generator is straightforward a function of the coordinates of the physical field, which could generate a response of any resolution. Then such a dimension reduction ratio would be infinite! To illustrate the resolution-independent advantage, a synthetic observation represented by 64 random sensors is used for parameter inversion. Figure 11 gives specific information about the synthetic observation and the posterior of the response under the given conditions. The number of sensors and position coordinates for this observation are inconsistent with those of the training set, but the generator does not need to be retrained to generate a response for the corresponding coordinates. Table 2 shows the specific results, and the accuracy and confidence interval are still within a reasonable range. The POD-based architecture still has lower errors at the high noise level. While the figures of the posterior distribution of the parameters are not presented here in this scenario.

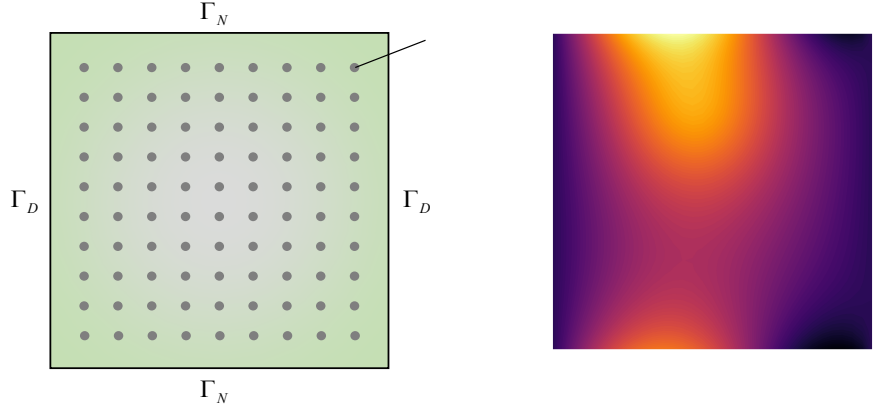

**Figure 3** Schematic diagram of Case 1. On the left is some information about the boundary and the sensor, and on the right is a high-resolution reference solution of the observation.

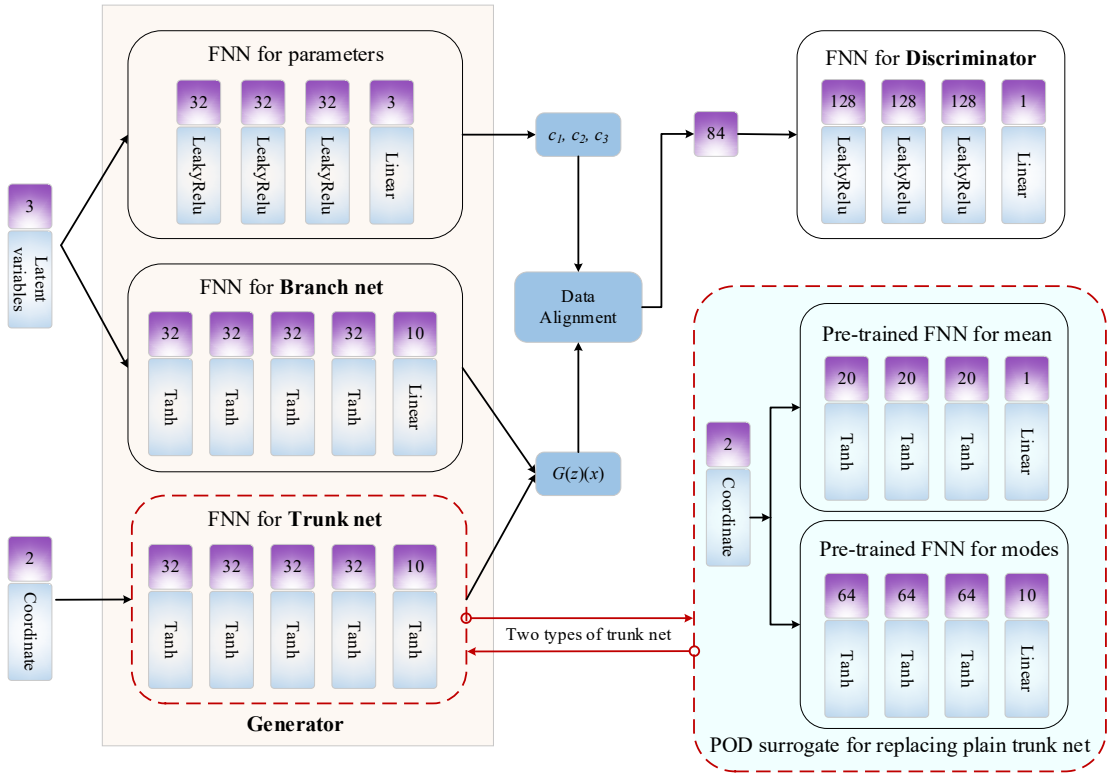

**Figure 4** The architecture of the neural networks in Case 1. The red dotted boxes represent the two available trunk net architectures.

Table 1 Identification results of the posterior distribution of $c_1$, $c_2$, and $c_3$. 81 sparse uniform sensors.

| method | parameters | noise level 0.01 | | | noise level 0.1 | | |
|---|---|---|---|---|---|---|---|
| | | mean | std | relative error | mean | std | relative error |
| vanilla | $c_1$ | 0.4526 | 0.0071 | 0.82% | 0.4400 | 0.0706 | 1.99% |
| | $c_2$ | 0.7433 | 0.0069 | 1.27% | 0.7546 | 0.0326 | 2.80% |
| | $c_3$ | 0.1081 | 0.0012 | 2.72% | 0.1080 | 0.0079 | 2.78% |
| POD-based | $c_1$ | 0.4551 | 0.0103 | 1.38% | 0.4461 | 0.0645 | 0.63% |
| | $c_2$ | 0.7436 | 0.0030 | 1.31% | 0.7516 | 0.0289 | 2.40% |
| | $c_3$ | 0.1138 | 0.0010 | 2.43% | 0.1122 | 0.0085 | 0.99% |

Table 2 Identification results of the posterior distribution of $c_1$, $c_2$, and $c_3$. 64 sparse random sensors.

| method | parameters | noise level 0.01 | | | noise level 0.1 | | |
|---|---|---|---|---|---|---|---|
| | | mean | std | relative error | mean | std | relative error |
| vanilla | $c_1$ | 0.4535 | 0.0077 | 1.03% | 0.4497 | 0.0697 | 1.73% |
| | $c_2$ | 0.7575 | 0.0070 | 3.21% | 0.7700 | 0.0266 | 4.91% |
| | $c_3$ | 0.1093 | 0.0014 | 1.60% | 0.1021 | 0.0094 | 8.04% |
| POD-based | $c_1$ | 0.4364 | 0.0097 | 2.80% | 0.4638 | 0.0663 | 0.63% |
| | $c_2$ | 0.7509 | 0.0043 | 2.30 % | 0.7675 | 0.0276 | 2.40% |
| | $c_3$ | 0.1107 | 0.0024 | 0.33% | 0.1089 | 0.0112 | 1.89 % |

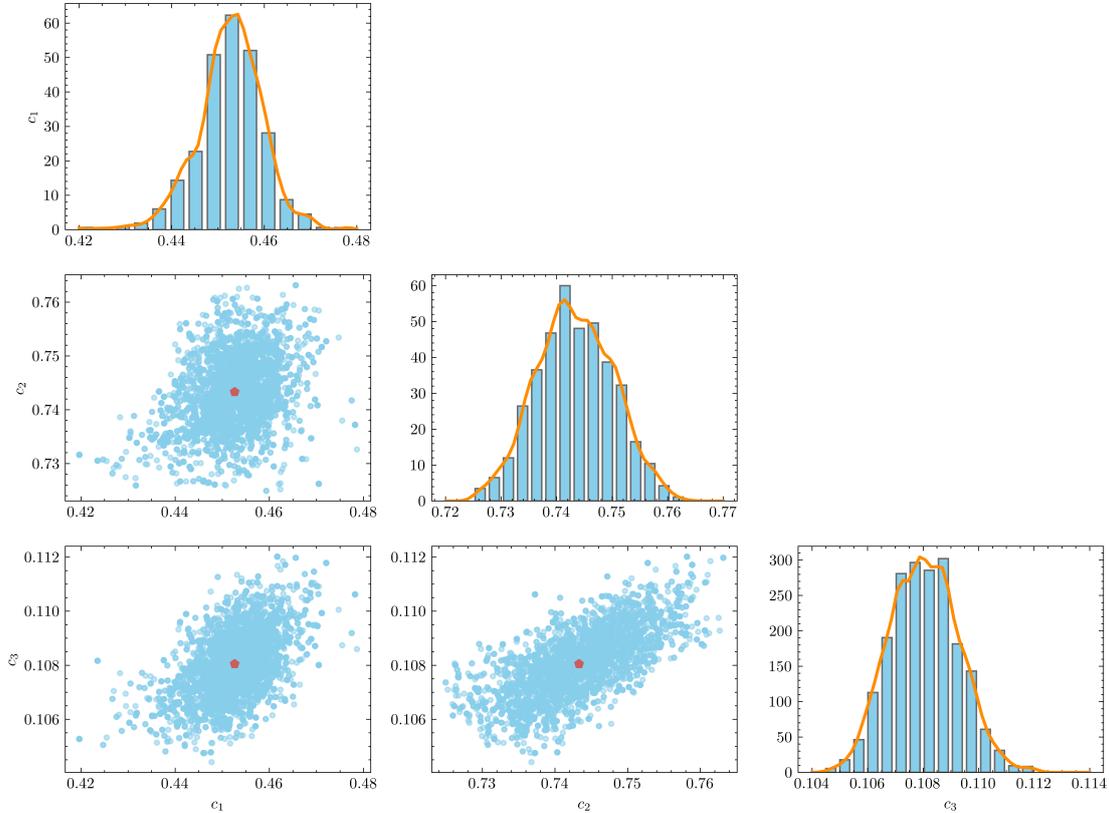

**Figure 5** Parameter identification results of Case 1, 81 uniform sensors, noise level 0.01, vanilla trunk net.

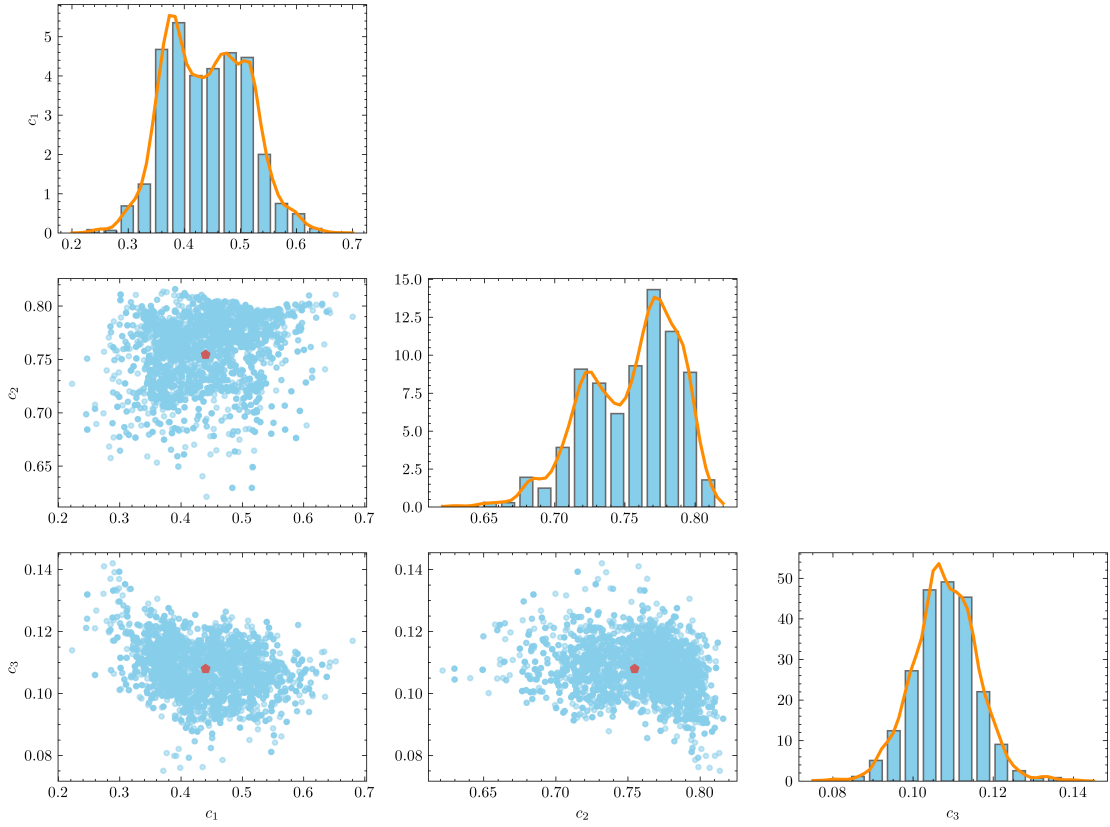

**Figure 6** Parameter identification results of Case 1, 81 uniform sensors, noise level 0.1, vanilla trunk net.

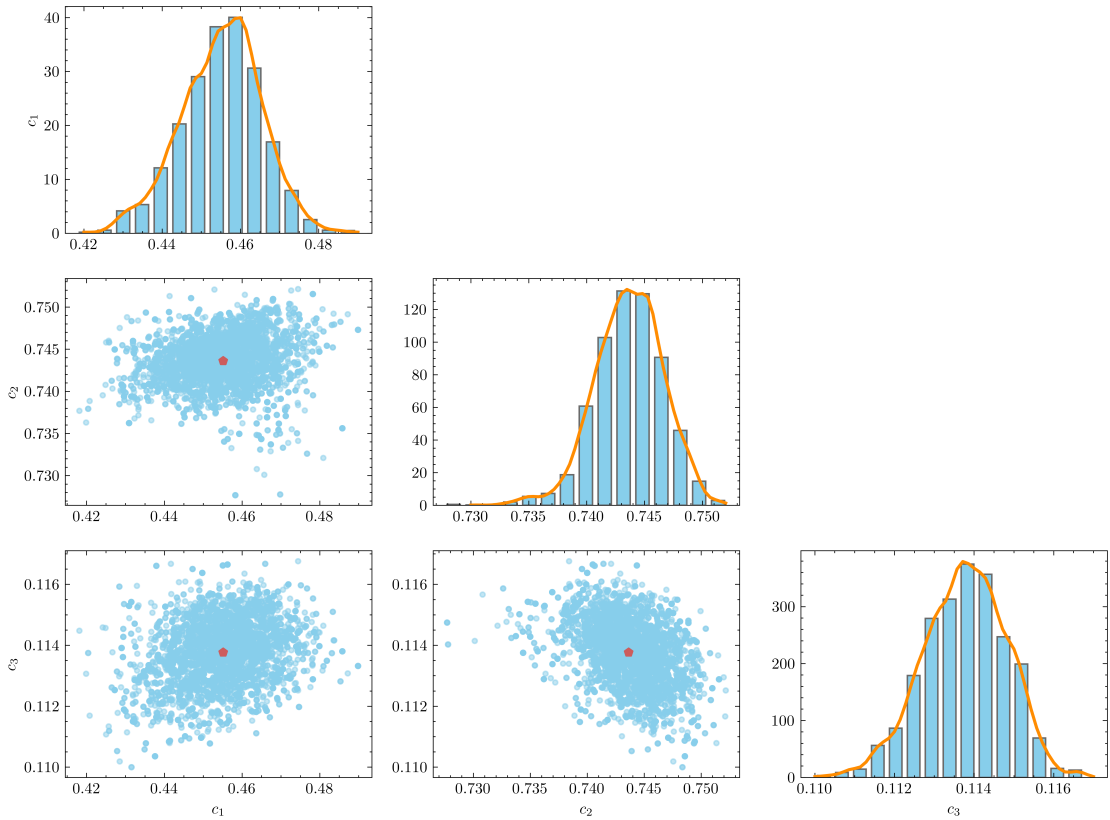

**Figure 7** Parameter identification results of Case 1, 81 uniform sensors, noise level 0.01, POD-based trunk net.

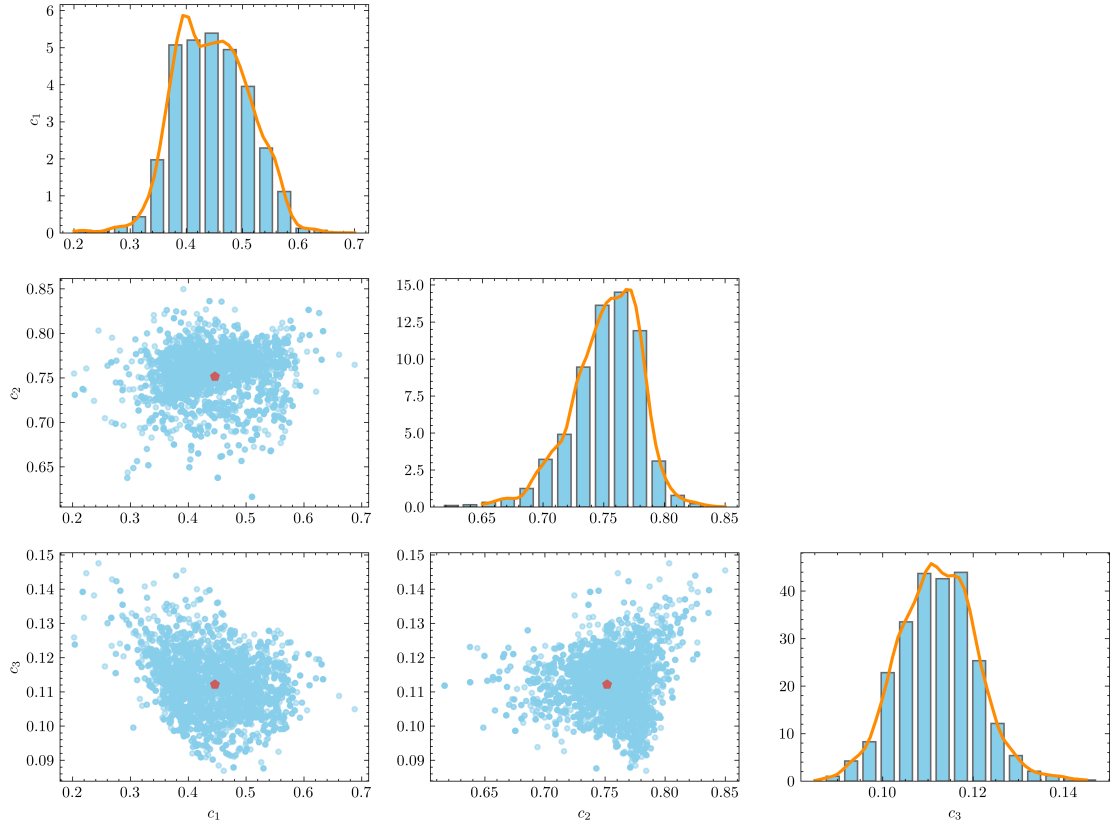

**Figure 8** Parameter identification results of Case 1, 81 uniform sensors, noise level 0.1, POD-based trunk net.

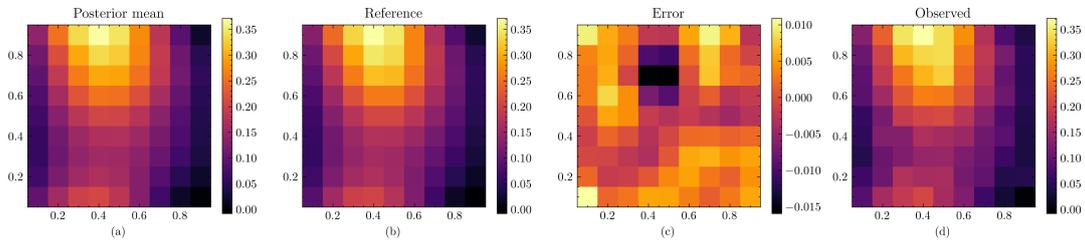

(A) noise level 0.01.

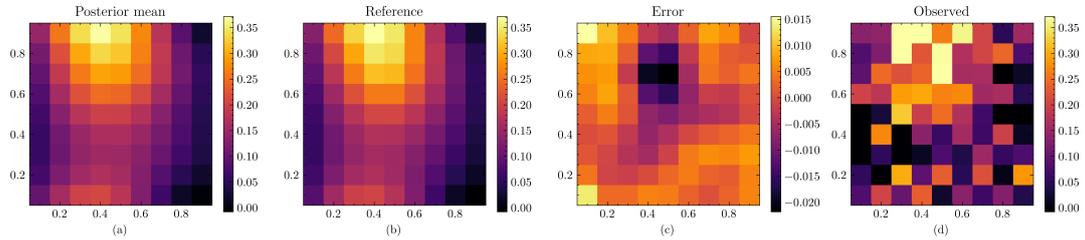

(B) noise level 0.1.

**Figure 9** Low-resolution response diagram composed of sparse sensors. The generated posterior mean is a "denoised version" of the observations. 81 uniform sensors, vanilla trunk net.

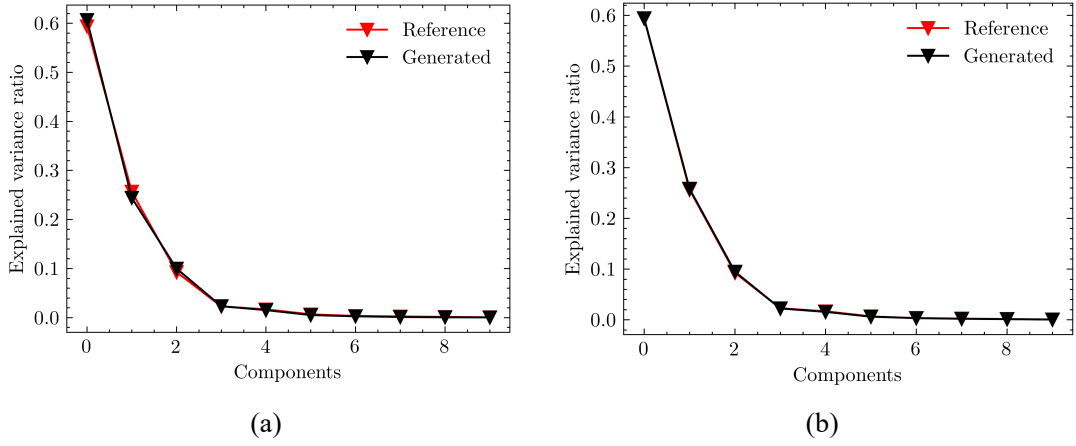

**Figure 10** The explained variance ratio of the reference responses and the generated responses of Case 1. (a) vanilla trunk net; (b) POD-based trunk net.

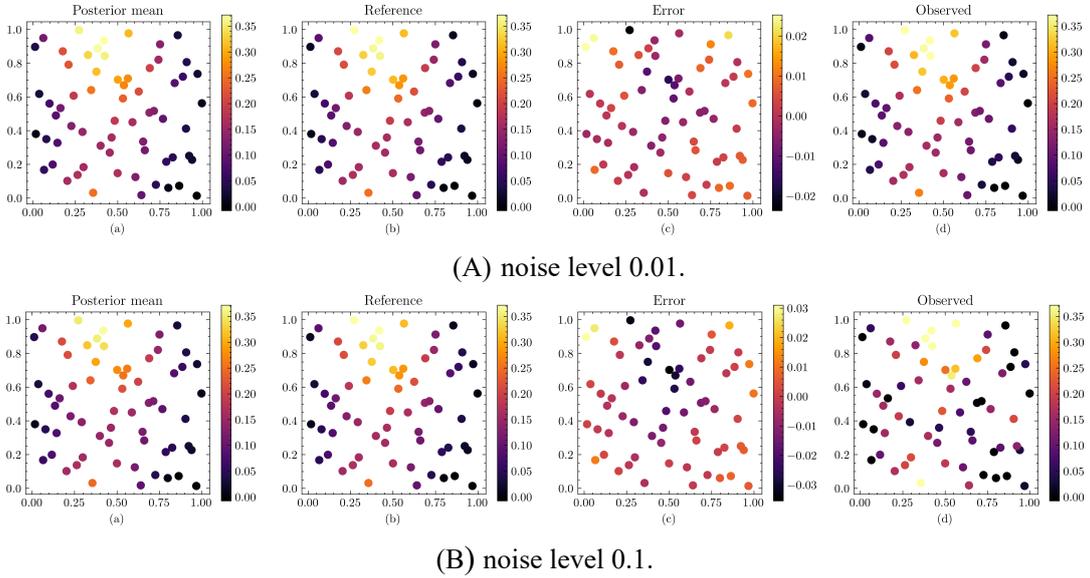

**Figure 11** Scatter diagram composed of sparse random sensors. The generated posterior mean is a "denoised version" of the observations. 64 random sensors, vanilla trunk net.

### 4.2 Case 2: elastic material inversion problem

The inverse problem considered in this section is a material identification problem in solid mechanics. The focus is placed on the estimation of spatially varying material properties under the assumption of dimensionless quantities. Specifically, the task here is to infer the elastic modulus $E(x, y)$ from the noisy observed deformation $u(x, y)$ under given conditions [18]. The problem description can be seen in Figure 12. It is a unit square domain $\Omega = [0,1] \times [0,1] \in \mathbb{R}^2$. In this case, the governing equation reads:

$$\nabla \cdot \boldsymbol{\sigma}(\boldsymbol{u}, E) = 0, \ x, y \in \Omega, \tag{22}$$

where $\sigma$ is the second-order stress tensor. Assuming that the problem satisfies the condition of isotropic and small strain theory [40]. Consequently, the stress can be described as a function of $E$ and $\boldsymbol{u}$ as follow:

$$\sigma = \frac{2E}{1+v}[\nabla \boldsymbol{u} + \nabla \boldsymbol{u}^T + \frac{2v}{1-v}tr(\nabla \boldsymbol{u})\mathbf{I}], \qquad (23)$$

where $v$ is the Poisson ratio and is supposed to be a constant of 0.3. $tr(\cdot)$ is the trace function and $\mathbf{I}$ is the identity tensor. Two types of BCs are applied with given conditions as follows:

$$\begin{aligned} u_x\big|_{x \in \Gamma_1} &= 0, \ u_y\big|_{x=0, y=0} = 0, \\ \sigma \cdot \mathbf{n}\big|_{\Gamma_3, \Gamma_4} &= 0, \\ \sigma \cdot \mathbf{n}\big|_{\Gamma_2} &= 1.5, \end{aligned} \qquad (24)$$

where $u_x$ and $u_y$ are the components of deformation in the $x$ and $y$ directions, respectively. $\mathbf{n}$ is the outward normal for the corresponding surface. The problem is described in Figure 11. The true elastic modulus is modeled by the Gaussian process as a random field:

$$\begin{cases} E(x,y) = 1 + 0.1\exp(g(x,y)), \\ g(x,y) \sim \mathcal{GP}(0, \exp(-\frac{1}{2}((x-x')^2 + (y-y')^2))). \end{cases} \qquad (25)$$

The samples of $E$ are not generated directly by this formulation, while a truncated KLE of $g(x,y)$, which retains 99.9% of the energy, is used and then the above transformation is applied. A data set of 990 samples of $E$ is used as a training set and the associated displacement was calculated by FEM. The domain is discretized to a $25 \times 25$ gird. Assuming that all the node displacements are the responses in this case, therefore the joint distribution of $E$ and $\boldsymbol{u}$ reaches a high resolution of 1,825. As before, an observation is taken from the test set, including the displacement in both directions and its elastic modulus needs to be inferred. The architecture of the OL-GAN in this case is shown in Figure 13. Notably, since the generator needs to generate a material field and two displacement fields, the operator network in this example is multi-output. There is more than one way to handle multiple function outputs and it is highly relevant to the problem [30]. In this study, we used 3 independent DeepONet and obtained reasonable results. The MSE of the pre-trained POD surrogate is less than 1e-7. The batch size is 330. The dimensionality of the latent space is 6. The results for this case are as follows.

The inversion results of the different scenarios are shown in Table 3. The performance of different generator architectures at different noise levels is compared. The resolution-independent inversion capability was also demonstrated by using 400 random sensors. The relative $l_2$ error is used to measure the error of the inversion result. It is found that in this example using a POD-based architecture can achieve relatively better results, more robust, with a minimum relative $l_2$ error of 8.6e-3 at the noise level of 0.01. In Figure 14 and Figure 15, specific observation information is given under different noises, i.e., displacement in the *x*- and *y*-directions, corresponding to (f) and (h). The posterior mean, absolute error, and sample standard deviation of the unknown parameters, namely the elastic modulus, are shown in the first row of these two figures. Meanwhile, the posterior corresponding to the observational information is also given. Figure 16 shows the result obtained by the vanilla version of the generator, where only the posterior information for the material field of interest is shown. From this set of comparison figures, it can be seen more intuitively that the result generated by the POD-based architecture is closer to the reference solution in the sense of a contour map since the error is lower. To illustrate that both architectures have learned the approximated overall distribution of the physical field data, Figure 17 is plotted by using fake samples and training samples. From the point of view of the explained variance ratio, the spectra of the generated samples are acceptable, but the displacement field is relatively more accurate because the correlation of the variables is simpler than that of the material field. The contour map of the elasticity modulus generated by the vanilla generator is not similar at a noise level of 0.01 with a relative $l_2$ error of 0.0386, indicating that the generator has not properly learned the local features in the domain. This may be due to the randomness of training GANs, or it may be due to the architecture of the generator, which will require more research in the future. Figure 18 and Figure 19 are the results of inversions using 400 random sensors, which also demonstrate the flexibility of the generator. Specifically, we can first use the coordinates corresponding to the random observation data to find the posterior distribution of the latent variables, and then use the latent variables and the coordinates of the regular grid to generate the material field and the displacement field. Finally, the contour map can be drawn. Only the results of the POD-based architecture are given here for illustration. In addition, we also observe that the network architecture based on operator learning can produce a smoother physical field, as detailed in Appendix B.

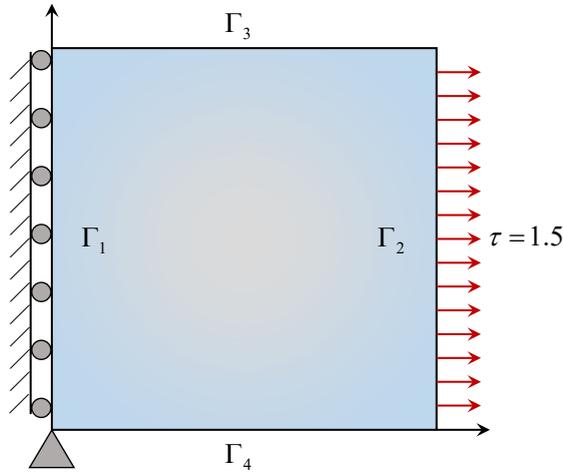

**Figure 12** Schematic diagram of Case 2.

**Table 3** Relative $l_2$ error of the posterior mean of the elasticity modulus field for different scenarios.

| noise level | 625 uniform sensors | | 400 random sensors | |
|---|---|---|---|---|
| | vanilla | POD-based | vanilla | POD-based |
| 0.01 | 0.0386 | 0.0086 | 0.0214 | 0.0096 |
| 0.1 | 0.0236 | 0.0103 | 0.0233 | 0.0173 |

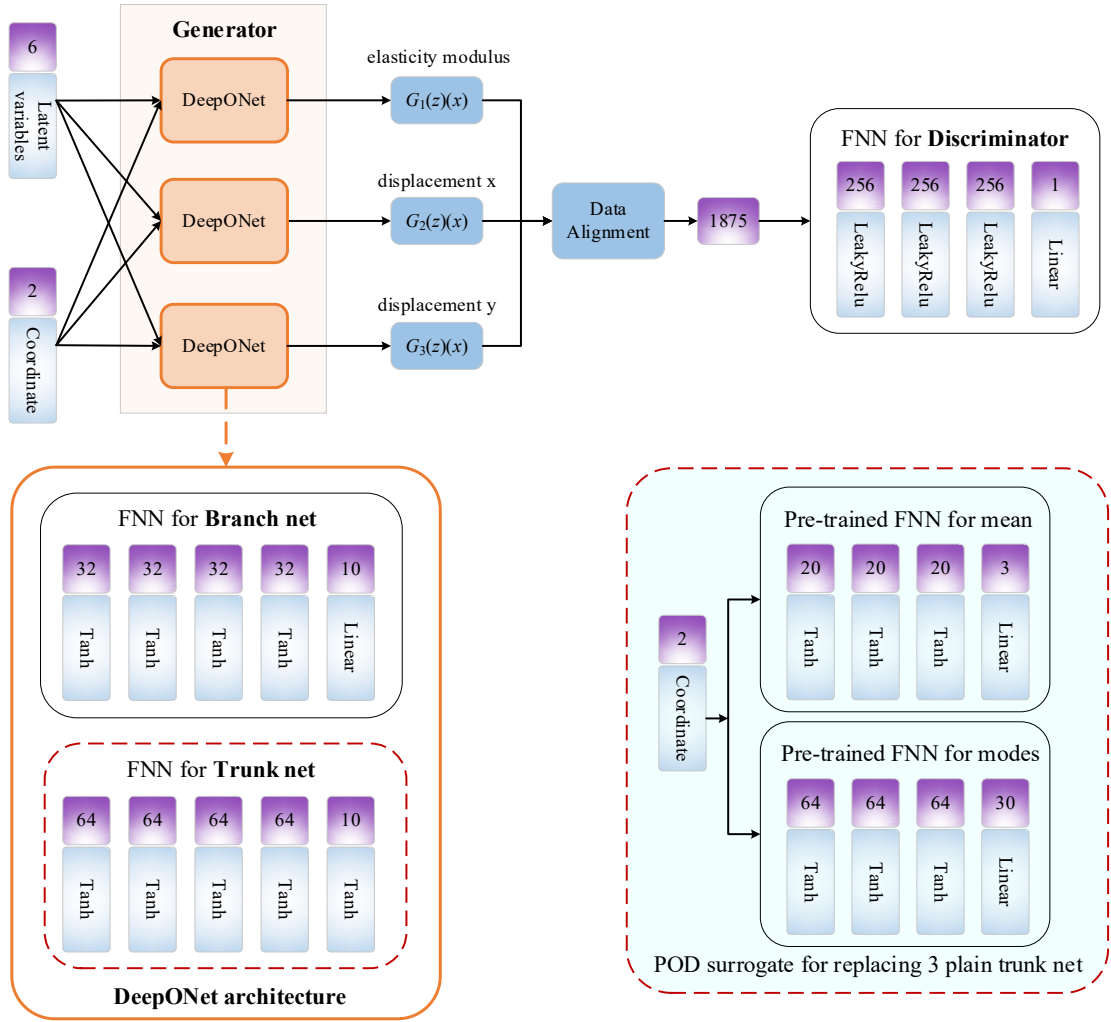

**Figure 13** The architecture of neural networks in Case 2. The red dotted boxes are two available trunk net architectures. In this example, the outputs of the POD surrogate include the mean and the first 10 modes of the 3 fields. Corresponding elasticity modulus $E$, $x$ displacement, $y$ displacement.

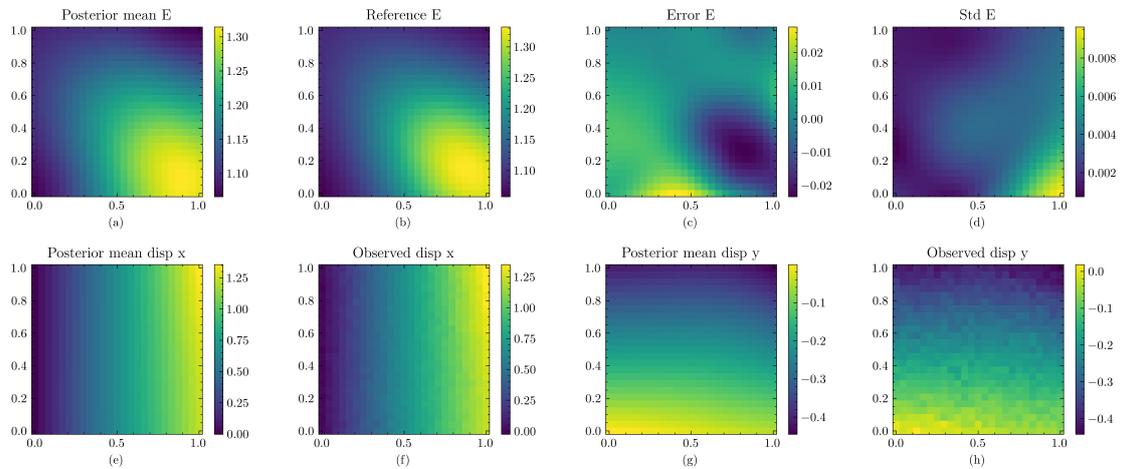

**Figure 14** Identification results and observed data of Case 2. 625 uniform sensors, noise level 0.01, POD-based trunk net.

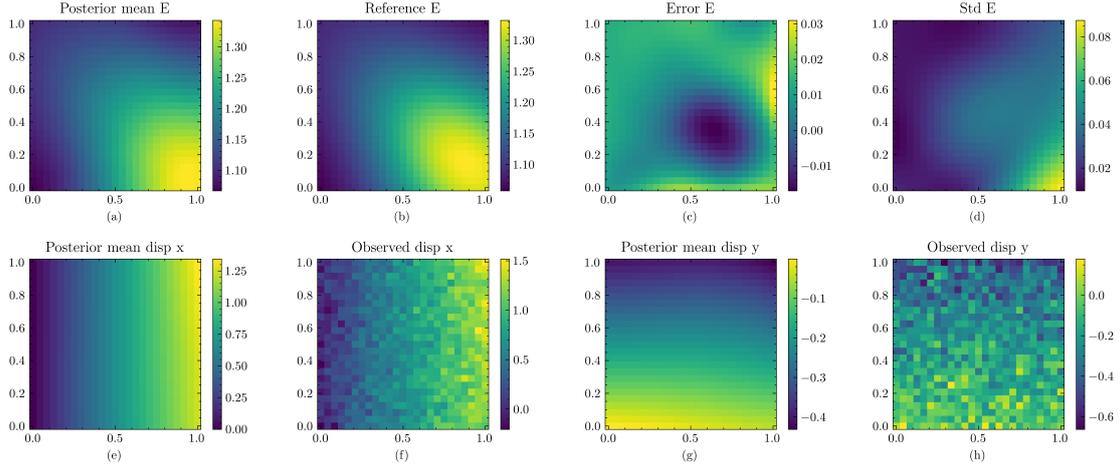

**Figure 15** Identification results and observed data of Case 2. 625 uniform sensors, noise level 0.1, POD-based trunk net.

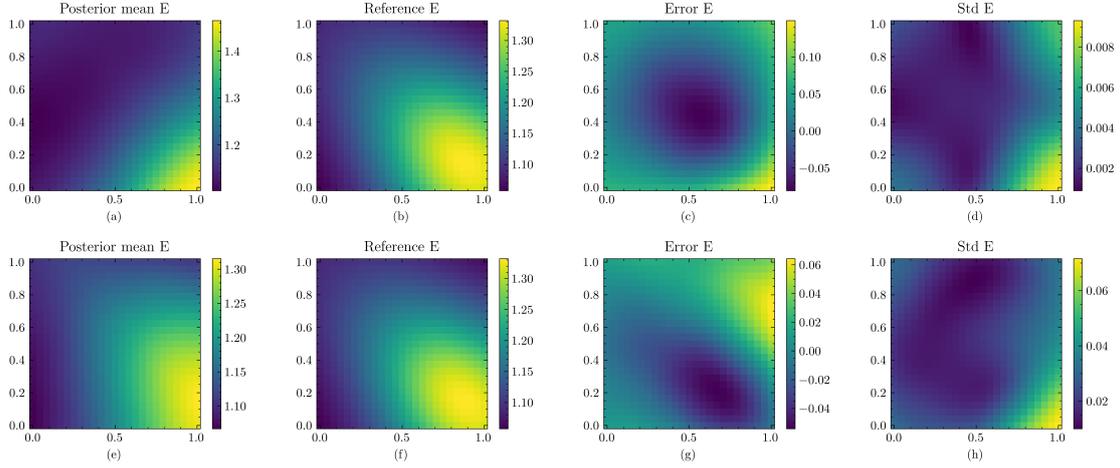

**Figure 16** Identification results of elasticity modulus of Case 2. 625 uniform sensors, vanilla trunk net. The noise of observations corresponding to the first row is 0.01, and the second row is 0.1.

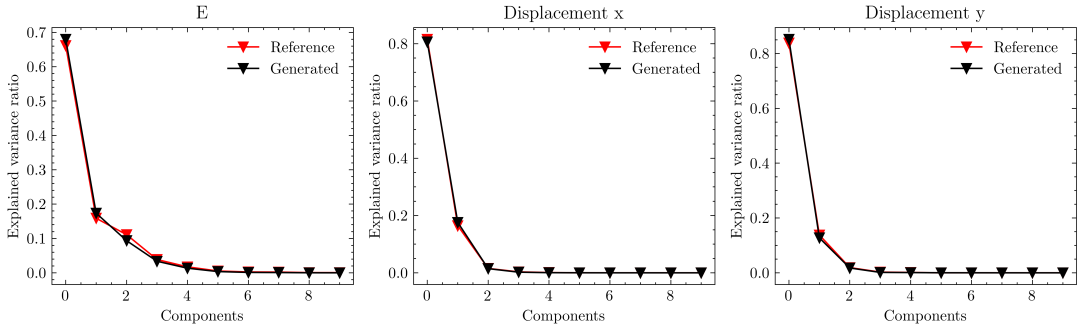

(a) vanilla trunk net.

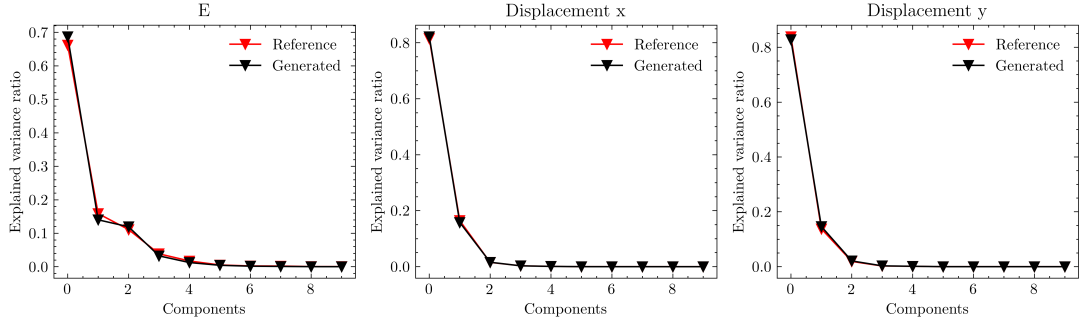

(b) POD-based trunk net

**Figure 17** The explained variance ratio of the reference and the generated samples of Case 2.

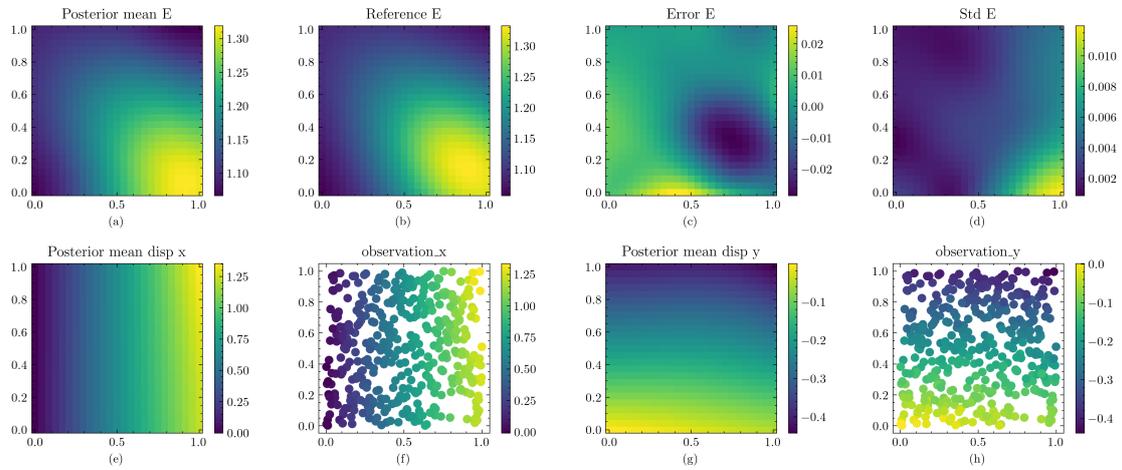

**Figure 18** Identification results and observed data of case 2. 400 random sensors, noise level 0.01, POD-based trunk net.

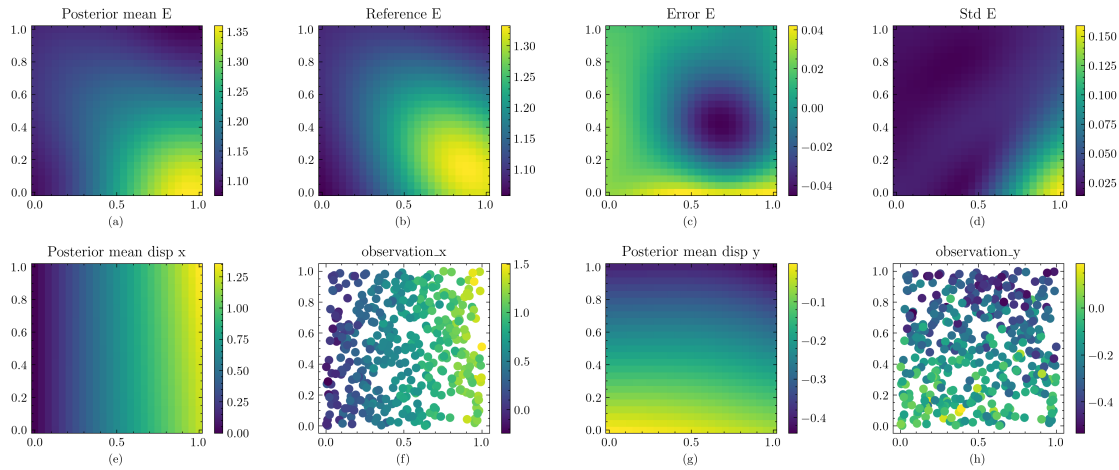

**Figure 19** Identification results and observed data of case 2. 400 random sensors, noise level 0.1, POD-based trunk net.

### 4.3 *Case 3: diffusion-reaction dynamic system*

Considering a diffusion-reaction system with a source term $u(x)$ which can be described by

$$\frac{\partial s}{\partial t} = D_s \frac{\partial s}{\partial x^2} + ks^2 + u(x), \quad x \in [0,1], t \in [0,1], \qquad (26)$$

where $D_s=0.01$ is the diffusion coefficient, and $k=0.01$ is the reaction rate. The initial conditions and boundary conditions are supposed to be zero. The task here is to infer the source term from the observed responses. The source term $u(x)$ is generated by sampling from a Gaussian process with zero mean and exponential quadratic covariance kernel with a length of 0.2. The source is represented by 100 uniformly distributed points over [0, 1]. A high-fidelity data set containing 2,000 samples is generated by solving the diffusion-reaction system using a second-order implicit finite difference method on a $100 \times 100$ grid. The training set is then obtained by interpolating on a $33 \times 33$ grid, meaning 33 uniformly distributed sensors with a time step 1/32. Hence, the resolution of joint distribution is 1,189. Consistent with previous treatment, assuming an observation $s(x, t)$ is taken from the test set, and two additive Gaussian noise levels of 0.1 and 0.2 are studied. The architecture of the OL-GAN in this case is shown in Figure 20. In this study, we used 2 independent DeepONet to generate parameters and responses. The MSE of the pre-trained POD surrogate is less than 1e-7. The batch size is 500. The dimensionality of the latent space is 6. The results for this case are as follows.

The inversion results of different scenarios are shown in Table 4. In this case, in addition to the relative $l_2$ error, $R^2$ is also used to evaluate the results. It can be seen from the results that the proposed method achieves the desired accuracy at different noise levels. Each result has a very high $R^2$ that is bigger than 0.96, the best one reaches 0.9899. The results of the POD-based architecture are slightly better than that of the vanilla one in this case. To illustrate that both architectures have learned the distribution of the parameters and responses, Figure 21 is plotted. And good or even better results can be obtained under random sparse measurements for both architectures. Figure 22 shows the specific details of the inversion results for different scenarios. The reference parameter is the red line, the random sensor location is given in the corresponding observation diagram. This example further demonstrates the resolution-independent advantage of the generator, as parameters, responses, and observations can all use different coordinates and, to a certain extent, do not affect the accuracy of the computation. The parameter uncertainty estimates are also reasonable. The main errors are concentrated at both endpoints, and the uncertainty is also greater near the two endpoints.

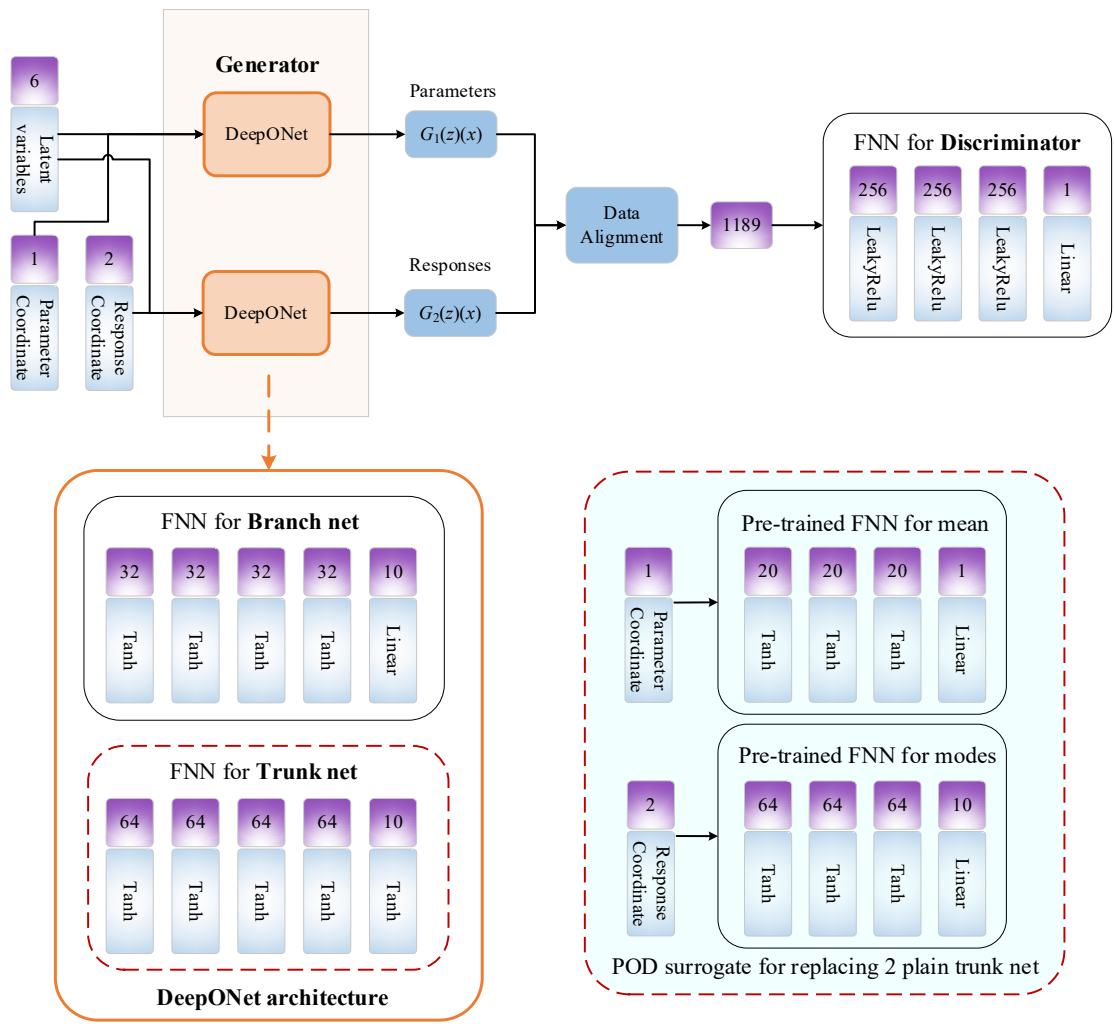

**Figure 20** The architecture of neural networks in Case 3. The red dotted boxes are two available trunk net architectures. In this example, the outputs of the POD surrogate include the mean and the first 10 modes of the 2 fields. Corresponding source terms and responses which is described by different discrete representations.

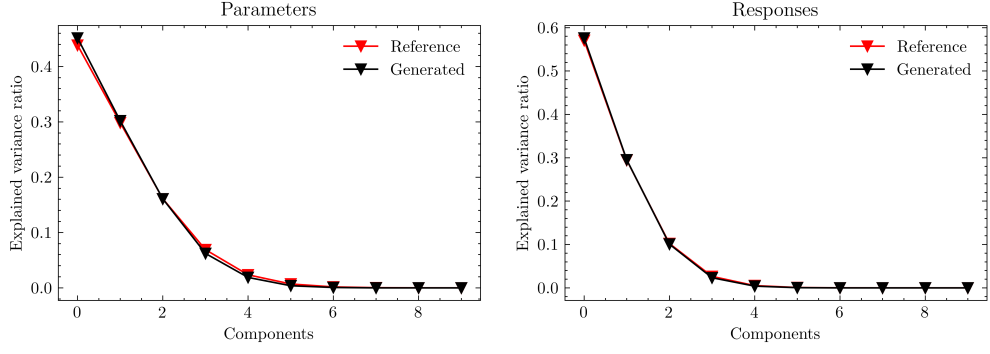

(a) vanilla trunk net.

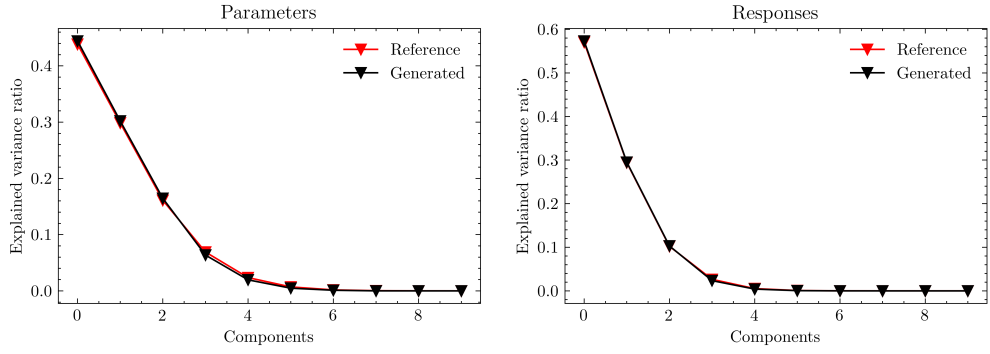

(b) POD-based trunk net

**Figure 21** The explained variance ratio of the reference and the generated samples of Case 3.

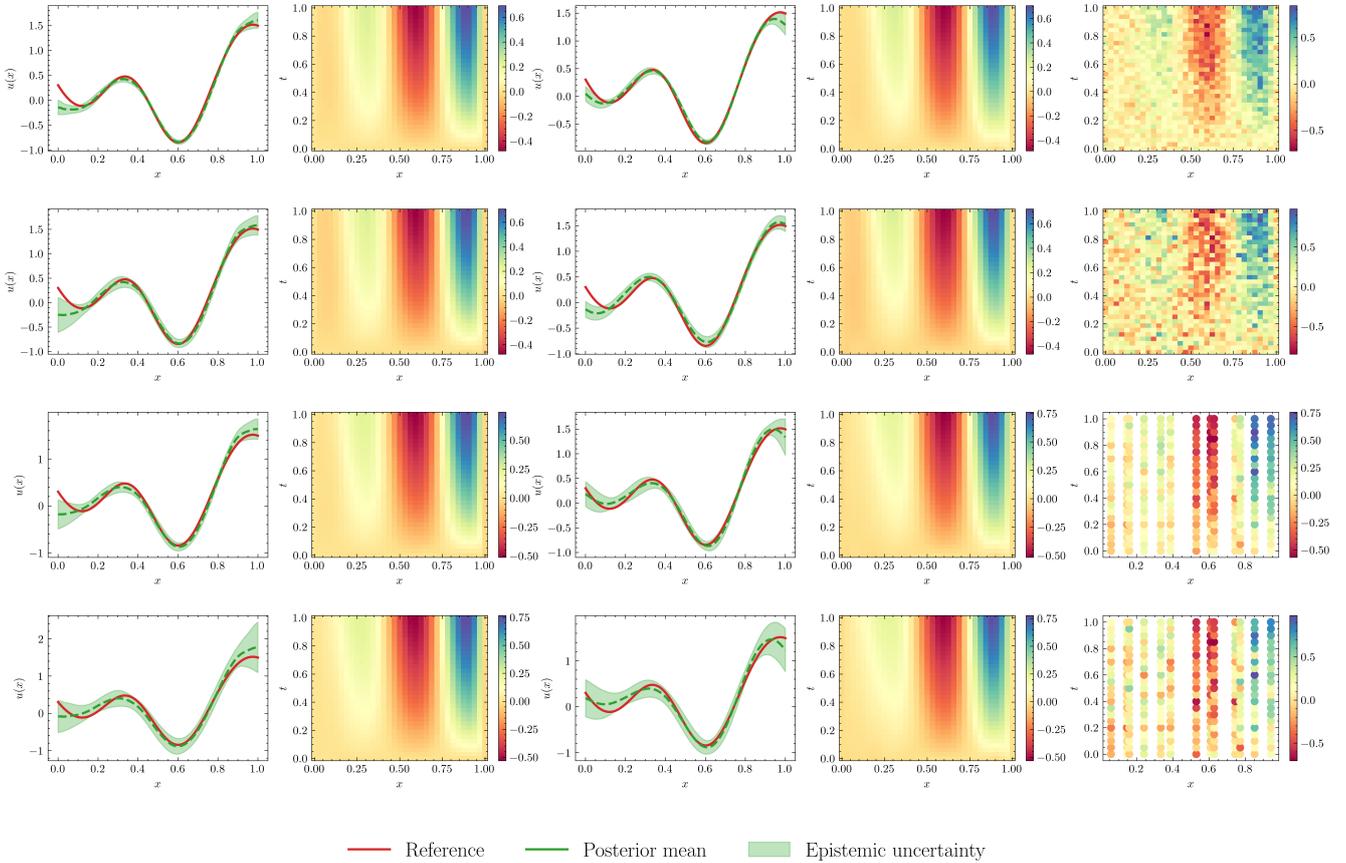

**Figure 22** Identification results and observed data of Case 3. The first two columns are the identified

parameters and the posterior mean of responses using the vanilla architecture, respectively. The 3rd and 4th columns are relevant results of the POD-based architecture. The last column is the observation data of different scenarios. The noise level of the 1st and 3rd rows is 0.1, the remained rows are 0.2.

Table 4 Relative $l_2$ error and $R^2$ of the posterior mean of the parameters for different scenarios.

| metrics | noise level | 33 uniform sensors (time step 1/32) | | 13 random sensors (time step 1/20) | |
| --- | --- | --- | --- | --- | --- |
| | | vanilla | POD-based | vanilla | POD-based |
| $R^2$ | 0.1 | 0.9795 | 0.9885 | 0.9706 | 0.9899 |
| | 0.2 | 0.9686 | 0.9736 | 0.9630 | 0.9780 |
| Relative $l_2$ error | 0.1 | 0.1377 | 0.1028 | 0.1648 | 0.0966 |
| | 0.2 | 0.1704 | 0.1561 | 0.1848 | 0.1426 |

# 5. Conclusion

In this study, Bayesian inverse problems are solved within the proposed GANs framework based on operator learning, namely OL-GANs. This method partially mitigates the challenges associated with prior selection and computational cost in Bayesian methods, and it is resolution-independent. The OL-GANs are trained by the joint distribution of parameters and responses of the physical system, with latent variables being shared by each component of the joint distribution. The latent variable distribution follows a simple isotropic standard Gaussian distribution, enabling efficient execution of the MCMC method in this latent space.

Compared with other methods, our approach allows for flexible sensor placement where both position and number of sensors used for synthesizing observation data can differ from those in the training set. Three classical examples in computational physics are used to demonstrate the effectiveness of the method, including a heat source parameter estimation, an elastic material field estimation, and a diffusion-reaction dynamic system source estimation problem. These three problems contain independent parameter estimation and function estimation, yielding reasonable uncertainty estimates under different noise levels. Notably, the POD-based generator yields relatively better results in these examples, probably because of its ability to fully exploit the information contained within the data. In addition, it is also discussed in the appendix that the joint distribution helps to reduce the dimension of the parameters and that the generated samples will be smoother based on operator learning.

While feasible, this method also has certain limitations and deserves improvement. For example, the operator learning introduces coordinate space, and the input of the network requires positional information that may not be present in existing training data. Although satisfactory results were achieved using simple FNNs to build all models, other techniques such as normalization, residual blocks, etc. could be explored to improve the reliability and accuracy.

## Appendix A. Additional results and discussion

In this section, an additional result from Case 1 is provided to illustrate the distribution transformation in a lower dimensional latent space, which involves the dimensionality reduction of independent parameters. Here, the 2 hidden variables follow an independent identically Gaussian distribution, while the 3 parameters are independent and uniformly distributed. The main difficulty is the need to effectively find a suitable transformation from manifold to independent parameters. The relevant results are computed and shown in Table A. 1, which implies that the model has still learned the joint distribution of parameters and responses. The spectra of the real and the fake samples are shown in Figure A. 1. At a noise level of 0.01, the accuracy of $c_1$ and $c_2$ is very high, but the accuracy of $c_3$ is slightly worse, with a relative error of 7.57%. At a noise level of 0.1, only $c_2$ has a lower relative error of 1.35%. Figure A. 2 and Figure A. 3 show the details of the posterior distribution of the parameters. As can be seen from the scatter plot, the pairwise joint distribution of the parameters is more complex. For a small noise of 0.01, the joint distribution of $c_1$ and $c_2$ is narrow. When the noise is 0.1, the pairwise distribution of the parameters becomes distorted and strange. This might be because low-dimensional manifolds fill higher-dimensional regions by folding and twisting themselves [19]. Therefore, the sample distribution in the high-dimensional space generated by the latent variable transformation becomes less smooth. This can also be analyzed in conjunction with the results in Section 4.1.

When the dimension of the latent variable is the same as that of the parameter, the accuracy is higher, and the pairwise distribution of the posterior parameters is significantly more reasonable, without distortion.

To illustrate the effect of dimensionality reduction, two other classical algorithms, PCA and AutoEncoder (AE), are compared. It should be noted that dimension reduction based on GANs is done by transforming latent variables, which are known low-dimensional simple distributions, into original parameters via a generator. The latent variables corresponding to the original parameters cannot be obtained directly from the generator. The dimensionality reduction of PCA and AE is non-generative. Original parameters are encoded by a trained encoder (principal component) into latent variables (principal component weights), which can also be decoded into original parameters by a decoder. It should be noted that it is not a completely fair comparison due to the use of joint distribution in the proposed method. Different methods are trained with 1000 samples. The hyperparameters and network architectures of the GAN-based method are the same as the vanilla version in Section 4.1 (except for the input dimension). Notably, the activation function of the generator is replaced by Tanh, since better results than LeakyReLU are observed in this scenario. Two principal components are selected for dimensionality reduction based on PCA. The encoder and decoder network architecture of AE is symmetric. The decoder is an FNN which is the same as the parameter generator in the Gan-based method. The hyperparameters of the Adam for training AE are a learning rate of 0.001, betas = (0.9, 0.999). The number of iterations is 10,000.

The results are shown in Figure A. 4. It can be seen that the proposed method can generate the most similar samples and effectively achieve dimension reduction. The generated samples of pure WGAN-GP are not good with the same training epochs and hyperparameters. We believe that the joint distribution of parameters and responses helps the networks to find more appropriate weights to transform the manifold in a certain sense. Although the design parameters are independent, the discrete representations of the responses are dependent. Therefore, it is relatively easy for the response generator to learn the pattern of transformation of the latent variable to the response. The use of the joint distribution means that the corresponding parameters are used as "labels", which helps the networks to simultaneously transform the latent variables into independent parameters. However, PCA and AE can only retain the information of the first 2 dimensions in the case of 2 latent variables. This is because the design space is a "cuboid" where the sample variance of $c_1$ and $c_2$ is larger, so PCA will retain information in these directions. Similar results were obtained with AE since the loss function is MSE, which essentially minimizes the empirical risk of reconstructing the parameters. Hence, the network tends to reconstruct the sample information in the directions of $c_1$ and $c_2$, and eventually falls into a local optimum. To quantify the quality of the generated (reconstructed) samples, the $W_1$ distance between the generated (reconstructed) sample and the real sample is calculated in each dimension:

$$W_1(P,Q) = \inf_{\pi \in \Gamma(P,Q)} \int_{\mathbb{R} \times \mathbb{R}} |X - Y| d\pi(X,Y) \tag{A. 1}$$

where $\Gamma(P,Q)$ is the set of all joint probability measures on $\mathbb{R} \times \mathbb{R}$ whose marginals are $P$ and $Q$ on the first and second factors, $X$ and $Y$, respectively. The results are shown in Table A. 2. It can be seen that in the $c_3$ direction, the $W_1$ distance of the samples generated by the proposed method is an order of magnitude smaller than that of the other methods. And the $W_1$ distance is more balanced in the 3 dimensions. Although the samples reconstructed by PCA and AE are very accurate in $c_1$ and $c_2$, the information in $c_3$ is almost completely lost.

Table A. 1 Identification results of the posterior distribution of $c_1$, $c_2$, and $c_3$.

| sparse uniform sensors (81) | noise level 0.01 | | | noise level 0.1 | | |
| --- | --- | --- | --- | --- | --- | --- |
| parameters | mean | std | relative error | mean | std | relative error |
| $c_1$ | 0.4553 | 0.0064 | 1.42% | 0.4145 | 0.0706 | 7.66 % |
| $c_2$ | 0.7332 | 0.0020 | 0.10% | 0.7439 | 0.0160 | 1.35% |
| $c_3$ | 0.1195 | 0.0006 | 7.57% | 0.1193 | 0.0063 | 7.36% |

Table A. 2 Comparison of $W_1$ distance of samples generated/reconstructed by different methods.

| direction | OL-GAN | pure WGAN-GP | PCA | AE |
| --- | --- | --- | --- | --- |
| $c_1$ | 0.00572 | 0.02424 | 0.00002 | 0.00127 |
| $c_2$ | 0.00575 | 0.03319 | 0.00006 | 0.00117 |
| $c_3$ | 0.00247 | 0.01104 | 0.02507 | 0.02438 |

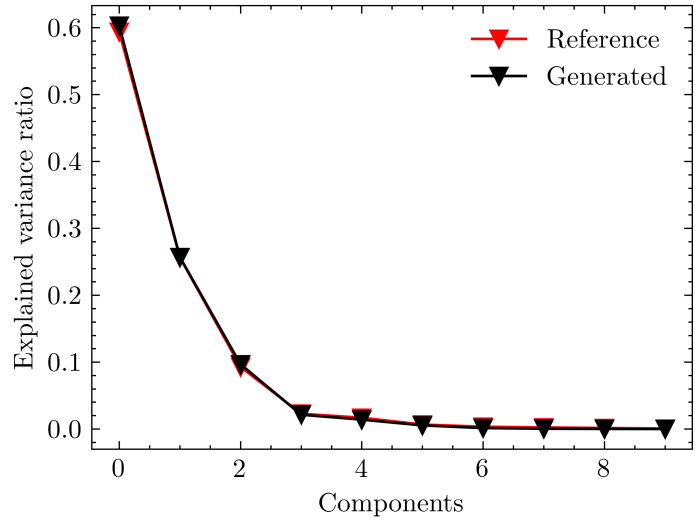

**Figure A. 1** The explained variance ratio of the reference responses and generated responses.

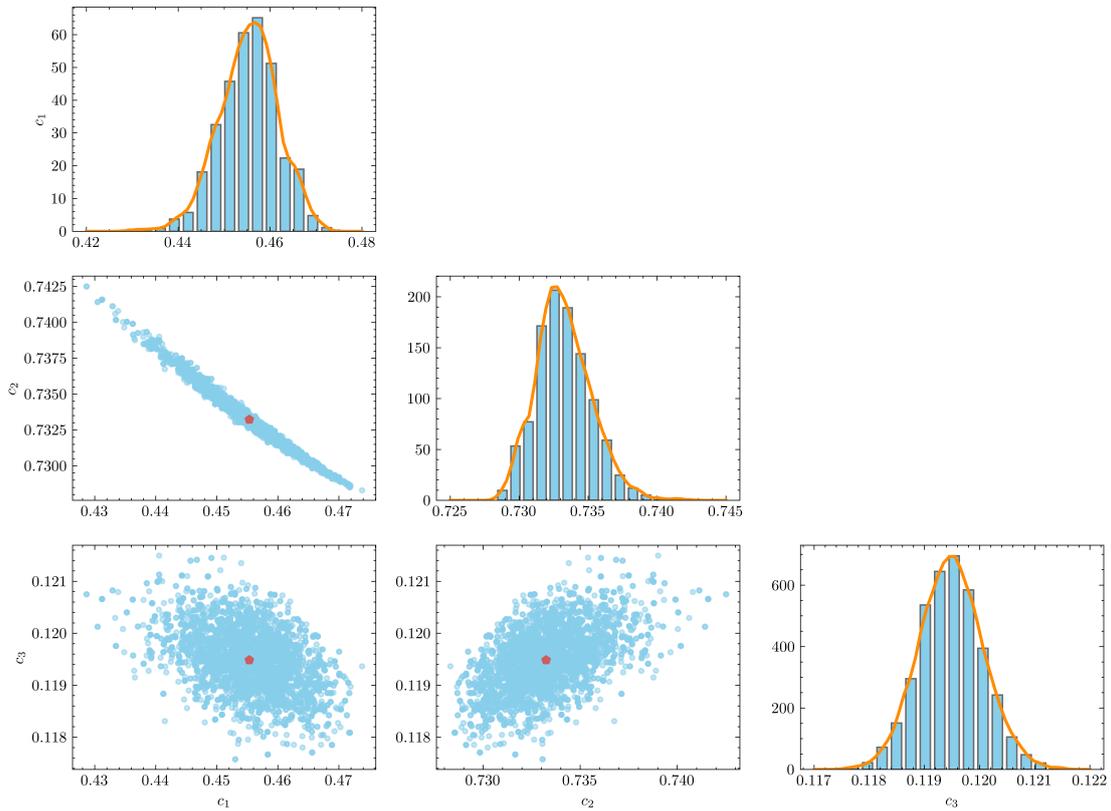

**Figure A. 2** Parameter identification results of case 1, 81 uniform sensors, noise level 0.01, vanilla trunk net, 2 latent variables.

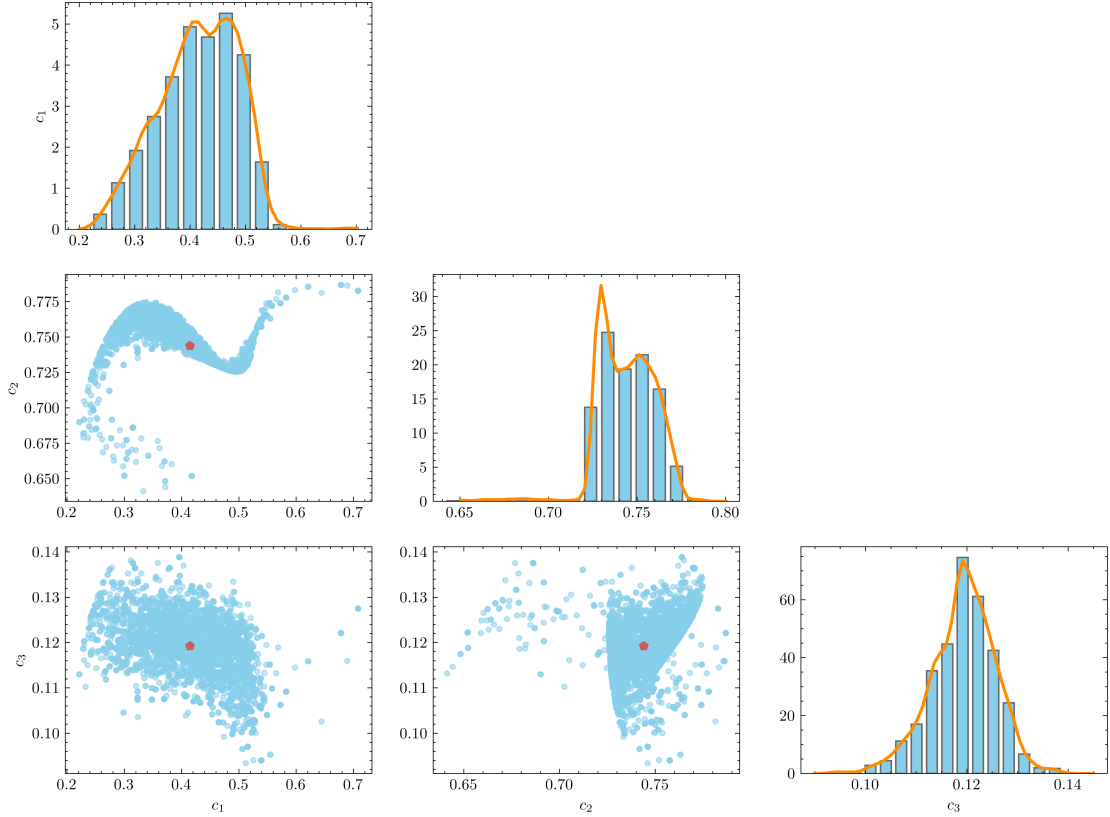

**Figure A. 3** Parameter identification results of case 1, 81 uniform sensors, noise level 0.1, vanilla trunk net, 2 latent variables.

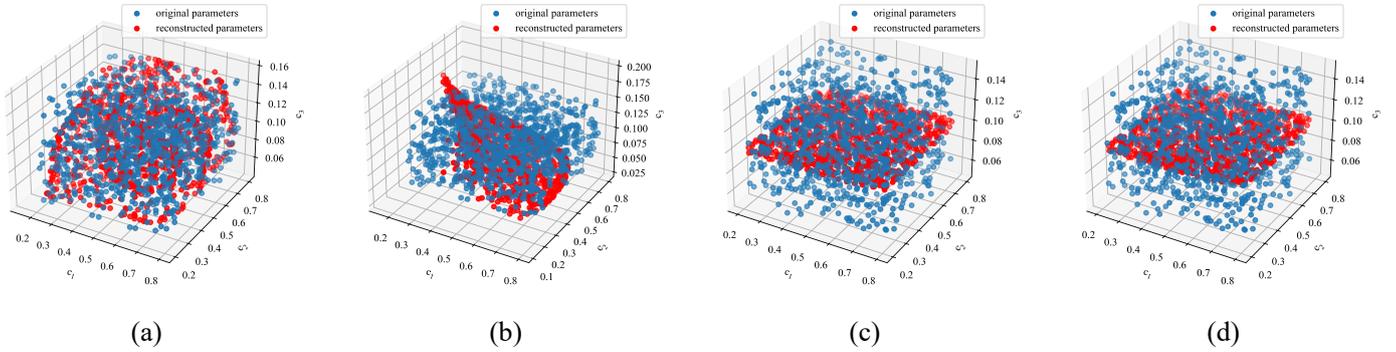

(a)           (b)           (c)           (d)

**Figure A. 4** Comparison of different methods for reconstructing. (a) Generated by the proposed method. The training samples are joint distributions (84 variables) and the network architecture is the same as the vanilla version in Section 4.1. (b) Generated by the pure WGAN-GP method. The training samples are only 3 independent uniformly distributed parameters and the architecture of the generator is only the " FNN for parameters " block of that in Section 4.1. (c) Reconstructed by PCA. (d) Reconstructed by the autoencoder.

# Appendix B. Smoothness

In this section, a result using a plain FNN as the generator is provided for comparison. Each physical field is generated by a network of four hidden layers. Each hidden layer has 32 neurons. The input dimension is 6 and the output dimension is 625. The capacity for such a network is comparable to those in Section 4.2. The other hyperparameters remain unchanged. As shown in the figure, there is a significant "noise" in the generated elastic modulus field with a relative $l_2$ error of 0.0231. While the results of previous generators based on operator learning are much smoother.

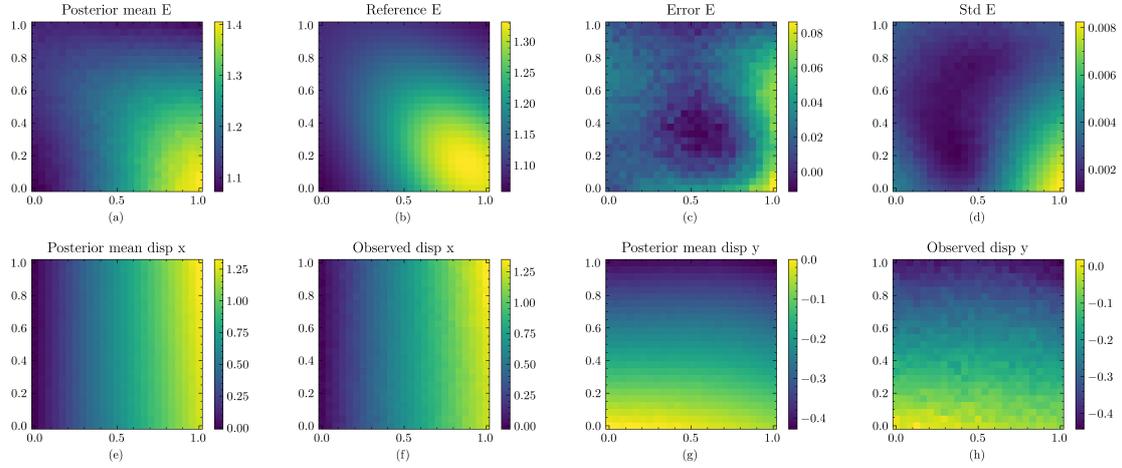

**Figure B. 1** Identification results and observed data of case 2. 625 uniform sensors, noise level 0.01, plain FNNs generators.

# Acknowledgments

This work has been supported by the National Key Research and Development Program of China (No. 2022YFB3303402), the Project of the National Natural Science Foundation of China (11972155), the Peacock Program for Overseas High-Level Talents Introduction of Shenzhen City (KQTD20200820113110016) and Project supported by Provincial Natural Science Foundation of Hunan(2020JJ4945).

# Declaration of interests

The authors declared that they have no conflicts of interest with this work. We declare that we do not have any commercial or associative interest that represents a conflict of interest in connection with this work.